\documentclass[12pt,notitlepage]{amsart}
\usepackage{amsmath, amssymb}     % [reqno] == equations numbered on right
\usepackage{amscd} 
\usepackage{tikz,tikz-cd}
\usepackage{import}
\usepackage{pdfpages}
\usepackage{xifthen}
\usepackage{transparent}
\usepackage{graphicx}
\usepackage{caption}
\usepackage{subcaption}
\usepackage[colorinlistoftodos]{todonotes}
\usepackage{hyperref}
\usepackage{enumitem}

\advance\voffset by -0.4 truein
\advance\hoffset by -0.4 truein
\newcommand{\Z}{\mathbb Z}

\newcommand{\cc}{\mathbb C} 
\newcommand{\qq}{\mathbb Q}
\newcommand{\p}{\mathbb P}

\newcommand{\g}{\mathfrak{V}}
\newcommand{\hh}{\mathfrak{W}}

\newcommand{\fx}{\mathfrak{X}}
\newcommand{\lp}{\mathbb{LP}}

\newcommand{\vp}{\varphi}
\newcommand{\ga}{\gamma}
\newcommand{\mo}{\mathcal{O}} 
\newcommand{\mv}{\mathcal{V}}

\renewcommand{\deg}{\mathrm{deg}\,}
\newcommand{\Ima}{\textnormal{Im}}

\DeclareMathOperator{\Id}{\textnormal{id}}

\theoremstyle{definition}
\newtheorem{theorem}[subsection]{Theorem}
\newtheorem{proposition}[subsection]{Proposition}
\newtheorem{Lem}[subsection]{Lemma}
\newtheorem{Cor}[subsection]{Corollary}

\newtheorem{Exa}[subsection]{Example}

\newtheorem{Def}[subsection]{Definition}
\font\smallrm=cmr8
\font\smallsc=cmcsc10
\font\smallsl=cmsl10

\begin{document}
\date{\today}
\author
[{\smallrm Eduardo Esteves, Renan Santos and Eduardo Vital}]
{Eduardo Esteves, Renan Santos and Eduardo Vital}
\title
[{\smallrm Quiver representations arising from degenerations of linear series, I
%simple bases for linked nets over $\Z^n$-quivers 
}]{Quiver representations arising from degenerations of linear series, I}
%characterising existence of simple basis; quiver representations akin to limit of linear series

\begin{abstract}
  We give a local characterization for when certain quiver representations in semisimple Abelian categories are semisimple, among them those arising from degenerations of linear series. 
  %the calculus of the multivariate Hilbert polynomial of the quiver Grassmannian associated to the representation. For these representations of quivers, we characterise the existence of simple bases with a  \emph{local} condition.
  This paper is the first of two, aimed to describe all the schematic limits of families of divisors associated to a given family of linear series on a one-dimensional family of projective varieties degenerating to a connected reduced projective scheme $X$ defined over any field, under the assumption that the total space of the family is regular along $X$, by means of certain quiver Grassmannians.
\end{abstract}

\thanks{E.~Esteves was supported by CNPq, Proc.~304623/2015-6 and FAPERJ, Proc.~E-26/202.992/2017, and R.~Santos and E.~Vital were supported by CAPES doctor grants at IMPA.}

\maketitle

\vspace{-.5cm}
\textbf{Keywords.} Linear Series \( \cdot \) Quivers \( \cdot \) Representations

\tableofcontents

% \section{Changes and improvements}
% \begin{enumerate}
%   \item Maybe we must cite \cite[Def. A.1.2]{osserman2014limit} where he define the \textit{Prelinked Grassmannian.} in our context a quiver representation in a general quiver.
%   \item Maybe it worth to cite \cite{10.1093/imrn/rnab241} here or in the Article II. Because maybe what we do here give another proof of some of their results.
% \end{enumerate}

\section{Introduction}\label{Intro}

This paper and its sequel \cite{Eduardo_Renan2021} aim to describe all the schematic limits of families of divisors associated to a given family of linear series on a family of projective varieties degenerating to a connected reduced projective scheme $X$ defined over any field $k$. We will describe these limits as points on a certain subscheme of the Hilbert scheme of $X$, and describe the subscheme using the quiver Grassmannian of pure dimension 1 of a certain quiver representation. 

%We study these quiver Grassmannians independently; we call them \emph{linked projective spaces}.

\emph{Linear series} are vector spaces of sections of a line bundle over a scheme over a field. A family of linear series on a family of smooth curves degenerating to a nodal curve gives rise to a \emph{limit linear series}, a term first coined by Eisenbud and Harris \cite{Eisenbud1986} for when the curve is of compact type: For each irreducible component $Y$ of the nodal curve, they considered the limit line bundle with degree 0 on all of the remaining components, and called its restriction and that of the limit sections to $Y$ an \emph{aspect} of the degeneration; the collection of all the aspects is the limit linear series. 

The theory found many applications, initially by Eisenbud and Harris themselves, \cite{EH2}, \cite{EH3}, \cite{EH4}, \cite{EH5} and \cite{EH6}, but later by many authors, as Farkas \cite{F2}, Cotterill \cite{Cot}, Tarasca \cite{Tar}, Farkas and Verra \cite{FV}, Pflueger \cite{Pfl}, and as recently as Knutsen and Lelli-Chiesa \cite{KLC}, to name a few.

It is thus no wonder that Eisenbud and Harris asked in \cite{EHBull}, p.~220, for a generalization of their theory to all nodal curves, writing that “...there is probably a small gold mine awaiting a general insight.” 

There have been attempts to generalize their theory, one even predating their work, by Ran \cite{Ran}, which was later independently recovered by the first author \cite{E0} and put to use in \cite{EM} and \cite{ES}. More recently, Tropical Geometry yielded another approach, starting with Baker and Norine \cite{BN} and Amini and Baker \cite{AB}, yielding results in \cite{CDPR}, \cite{JP1} and \cite{JP2}

Also, Osserman \cite{Osserman_2006} argued that for functorial reasons and for allowing work in positive characteristic, more ``aspects" should be considered. For curves of compact type with two components, he considered the collection of all limit line bundles of non-negative degree on each component, together with the limit sections, and called it a ``limit linear series." A decade later, Osserman \cite{Osserman_2019} studied general nodal curves, where a certain quiver was introduced; see loc.~cit., Notation~2.11, p.~62. It is essentially the quiver we use here, at least for when $X$ is a nodal curve, except that we consider only the trivial chain structure; see loc.~cit., Def.~2.2, p.~60. 

In contrast, our approach works for higher dimensional varieties. This may sound surprising, as there have been few works so far in this generality, but our approach works naturally.

We argue that a family of degenerating linear series on a family of projective varieties degenerating to a connected reduced projective scheme $X$ over a field $k$ gives rise to a certain quiver $Q$, a certain representation $\mathfrak L$ of $Q$ in the category of line bundles over $X$ and a certain subrepresentation of the induced representation $H^0(X,\mathfrak L)$ of $Q$ in the category of vector spaces over $k$ obtained by taking global sections in $\mathfrak L$; see Sections~\ref{lls_qui_rep_sb}~and~\ref{lnets}. The scheme $X$ needs not be a curve, $k$ needs not be algebraically closed nor of characteristic zero. These representations are barely studied in \cite{Osserman_2019}. We study them and prove:
\begin{enumerate}
    \item $Q$ is a \emph{$\mathbb Z^n$-quiver},
    \item $\mathfrak L$ is an \emph{exact} \emph{maximal} \emph{linked net}, 
    \item $\mathfrak V$ is a \emph{pure} \emph{exact} \emph{finitely generated} \emph{linked net}. 
\end{enumerate}
These statements are made in 
Section~\ref{degeneration_linear_series} and proved in Section~\ref{lnets}. The definition of a $\mathbb Z^n$-quiver appears in Section~\ref{Znq}, and the remaining definitions in Section~\ref{lnets}. The definitions, except for maximality, apply to any representation of $Q$ in a $k$-linear Abelian category, and all the arguments we use in this paper apply in this more general setting. 

Thus, having abstracted the properties the representations arising from degenerating linear series satisfy, we study here representations $\g$ of $\mathbb Z^n$-quivers $Q$ in a general $k$-linear Abelian category $\mathcal A$. The results we obtain here are essential though in \cite{Eduardo_Renan2021}, where we will focus more on representations by vector spaces and line bundles to obtain our goal of describing schematic limits of divisors. We refer to the introduction to \cite{Eduardo_Renan2021} for more motivation and for the road ahead. 

We will now proceed to explain what types of quivers and representations we study and what our main result in the current paper is.

A \emph{$\mathbb Z^n$-quiver} $Q$ is a quiver endowed with a nontrivial partition $A_0,\dots,A_m$ of its set of arrows $Q_1$ satisfying certain properties. Each $A_i$ is called an \emph{arrow type} and $a\in A_i$ is said to be of type $A_i$. The types satisfy certain properties: First, one can choose any finite sequence of arrow types and any vertex of $Q$ and obtain a unique path leaving that vertex whose sequence of arrows are of the types prescribed. Furthermore, if we reorder the sequence of arrow types, the path changes but the final vertex is the same. Finally, each vertex is connected to each other by an  \emph{admissible path}, that is, a path whose arrows are not of all types, and a path connects a vertex to itself if and only if it contains the same number of arrows of each type; see Section~\ref{Znq}.

A representation $\g$ of $Q$ in a $k$-linear Abelian category $\mathcal A$ is said to be a \emph{weakly linked net} over $Q$ if the composition $\vp^\g_\gamma$ of the maps associated to an admissible path $\gamma$ depends only on the vertices the path connects, up to homothety, and that associated to a non-admissible path is zero. We are interested in \emph{linked nets}, where the kernels of the two maps associated to two admissible paths having no arrow type in common intersect trivially. This is an ``open property." Another important open property is exactness. The representation $\g$ is \emph{exact} if the image of $\vp^\g_\gamma$ is equal to the kernel of $\vp^\g_\mu$ for any admissible paths $\gamma$ and $\mu$ whose concatenation $\mu\gamma$ is defined and is a nontrivial minimal circuit; see Section~\ref{lnets}.

%A representation having all the above mentioned properties is very special, and though many our results apply to less qualified representations, we focus on those. 

It is likely already apparent that the $\mathbb Z^n$-quiver $Q$ is infinite. We restrict though to weakly linked nets $\g$ that are \emph{finitely generated}, that is, for which there is a finite set of vertices $H\subset Q_0$ such that for each vertex $v$, the sum $\oplus \vp^\g_\gamma$ over all admissible paths $\gamma$ connecting $H$ to $v$ is an epimorphism. 

We work over the infinite quiver, but if one wanted to restrict $\g$ to $H$, one had better choose $H$ equal to its \emph{hull} $P(H)\subset Q_0$; see Definition~\ref{hull}. If $H$ is finite, so is $P(H)$ and the operation $\mathcal P(\cdot)$ is idempotent; see Proposition~\ref{HPH}. If $H$ generates $\g$, for each vertex $v$, there is a path $\gamma$ connecting $P(H)$ to $v$ such that $\vp^\g_\gamma$ is an epimorphism; see Proposition~\ref{simpfin}. 
We say $\g$ is \emph{$1$-generated} by $P(H)$. 

A vertex is equal to its hull. More generally, a \emph{polygon}, that is, a collection of vertices on a nontrivial minimal circuit of $Q$, is equal to its hull; see Definition~\ref{defpolygon} and Propositions~\ref{orientedmgon}~and~\ref{hullpolygon}.

A possibly weaker condition than being finitely generated is being locally finite: a weakly linked net $\g$ is \emph{locally finite} if for each vertex $v$ there is an integer $\ell$ such that $\vp^\g_\gamma=0$ for every path $\gamma$ arriving at $v$ of length greater than $\ell$; see Definition~\ref{deflnvsfs} and Proposition~\ref{noseq}. We classify locally finite exact linked nets of simple objects; they are those $1$-generated by a single vertex; see Theorem~\ref{simple1}. 

Put it differently, locally finite exact linked nets of simple objects are themselves the simple objects in the additive category $\mathcal C_Q(\mathcal A)$ of locally finite exact linked nets over $Q$ of objects in $\mathcal A$. It is a natural question to ask what the semisimple objects in $\mathcal C_Q(\mathcal A)$ are. We characterize them when $\mathcal A$ is semisimple, for example, when $\mathcal A$ is the category of finite-dimensional vector spaces. This is the main result of the paper, our Theorem~\ref{cat}..

%A standard question in Representation Theory can then be posed in our context: Is a pure exact finitely generated linked net of vector spaces \emph{completely decomposable}, that is, a direct sum of pure exact finitely generated linked subnets of dimension 1? Equivalently, do pure exact finitely generated linked nets admit simple bases? see \textcolor{orange}{ Definition~\ref{def_SB} and Theorem~\ref{simple2}}. This question is the focus of the current paper.

%A linked net has a simple basis if the representation decomposes as a direct sum of irreducible subrepresentations of pure dimension 1; see Definition~\ref{def_SB}. Linked nets arising as above for a curve $X$ with two components have simple bases. 

The characterization is not difficult if $n=1$; see Corollary~\ref{to_be_stated}: all objects of $\mathcal C_Q(\mathcal A)$ are semisimple! Not surprisingly, this result appeared already in the literature in various forms: Eisenbud and Harris asserted it in \cite[Lem.~2.3]{Eisenbud1986}, when they show the existence of an adapted basis; Osserman used it in his proofs of Lem.~A.9 and Lem.~A.12 in \cite{Osserman_2006}, and later in \cite{Osserman_2019}; it appeared more explicitly in Esteves and Osserman \cite[Lem.~2.3]{Esteves2013}, as well as more subtly in the proofs by Mu\~noz of \cite[Prop.~4.1]{G_munoz} and \cite[Prop.~3.6]{G_munoz_2020}. Osserman himself identified it as an important property, and essentially asked the first author in a private communication whether the result holds for higher $n$, at least for the linked nets arising from degenerations of linear series. 

The answer to the question is actually no, even for $\mathbb Z^2$-quivers. In Example~\ref{exac_no_SB} we present a finitely generated exact linked net of vector spaces of dimension $2$ over a $\mathbb Z^2$-quiver that is not semisimple. Furthermore, we show that the linked net arises from a degeneration of linear series along a pencil of cubics! 

There are though important cases where locally finite exact linked nets are semisimple, for instance, when they are generated by polygons; see Proposition~\ref{simpletriangle}. 
We will explain below why such linked nets are important. 

For the results in Section~\ref{applications}, we use our main result, Theorem~\ref{cat}, a characterization of the locally finite exact linked nets of objects in a $k$-linear semisimple Abelian category which are semisimple. It is a local characterization, in the sense that it is checked at each vertex of the quiver by comparing the kernels of the maps associated to the paths of length at most $n$ leaving the vertex; see Definition~\ref{defintprop}. Furthermore, it is a computable characterization: If the linked net is finitely generated, it is enough to check the property at the vertices of the finite set generating it; see Proposition~\ref{intH}.  We use the theorem in  Example~\ref{exac_no_SB} and Proposition~\ref{simpletriangle}.

%The known fact that linked nets arising from degenerations to a curve with two components have simple bases can be easily derived from it. We use this characterization to give an example of the non-existence of simple bases in the three-component case, our Example~\ref{exac_no_SB}. The example answers a private question made by Osserman to the first author in 2017.

%By doing away with the line bundles we focus on the linear aspects of a limit linear series. One of the goals of this paper is to present the basis of the theory of the special quivers and representations that arise. The fundamentals we study here will be applied in our second paper \cite{Eduardo_Renan2021}, in the study of \emph{linked projective spaces}. 

In the sequel \cite{Eduardo_Renan2021} we study the \emph{linked projective space} $\lp(\g)$ associated to a pure exact finitely generated linked net $\g$ of vector spaces over a $\mathbb Z^n$-quiver. The name is borrowed from work by Osserman on linked Grassmannians \cite{Osserman_2006}. In fact, $\lp(\g)$ is simply the quiver Grassmannian of subrepresentations of $\g$ of pure dimension 1. We describe its irreducible components, which are rational of the same dimension, and characterize when a number of them intersect, using polygons. We prove that $\lp(\g)$ is a local complete intersection and reduced.
%, arguing that $\lp(\g)$ is everywhere locally isomorphic to an open subscheme of $\lp(\g')$ for a linked net $\g'$ generated by a polygon, and proving directly that $\lp(\g')$ is Cohen--Macaulay. 
We prove that $\lp(\g)$ can be seen as a degeneration of the small diagonal in a product of projective spaces of equal dimension if $\g$ is smoothable, which is the case when $\g$ arises from a degeneration of linear series. 

Finally, in \cite{Eduardo_Renan2021} we give our own version of what a ``limit linear series" is, by defining a \emph{linked net of linear series} on a connected reduced projective scheme $X$ as the data $\mathfrak g=(Q,\mathfrak L,\g)$ of a $\mathbb Z^n$-quiver $Q$, a maximal linked net $\mathfrak L$ of line bundles over $X$ and a pure finitely generated linked subnet $\g$ of the linked net of spaces of global sections of the line bundles in $\mathfrak L$. If $\g$ and $\mathfrak L$ are exact, we expand on work by Rocha \cite{Rocha2019} for $\mathbb Z^1$-quivers to show that there is a natural morphism $\lp(\g)\to\text{Hilb}_X$ whose image is the set of schematic limits of divisors if $\mathfrak g$ arises from a degeneration.

This paper is organised as follows. In Section~\ref{Znq} we define $\mathbb Z^n$-quivers. We show in Section~\ref{lls_qui_rep_sb} how $\mathbb Z^n$-quivers and certain representations of them arise from degenerations of linear series, motivating the definition of linked nets in Section~\ref{lnets}.  Basic properties of $\mathbb Z^n$-quivers are presented in Section~\ref{generalquiver} and of linked nets in Section~\ref{generalnets}.
%We show in Theorem~\ref{Znquivers} that linked nets can be very different one from the other, but all $\Z^n$-quivers are isomorphic.  
%Implying that when the nodal curve is of compact type or not the associated {$\Z^n$-quiver} is the same, up to isomorphism. In addition, we give visual examples of a {$\Z^1$-quiver} and a $\Z^2$-quiver.
In Section~\ref{sec_SB} we define binary linked nets, of which a crucial example are the simple linked nets, whose associated objects are simple. We prove that simple locally finite exact linked nets are $1$-generated by a vertex. In Section~\ref{intprop} we introduce the \emph{intersection property}, which is a local condition for a linked net at a vertex of the quiver. 
%We also show some results in preparation to prove our main theorem, as Proposition~\ref{thm_quo_inter_pro}. 
In Section~\ref{Classification} we prove Theorem~\ref{thm_geneal_case} and Theorem~\ref{cat}: A locally finite exact linked net of objects in a $k$-linear semisimple Abelian category is semisimple it it satisfies everywhere the intersection property, which we apply in Section~\ref{applications} to show that an exact linked net generated by a polygon is semisimple, our Proposition~\ref{simpletriangle}, and to explain an example of an exact linked net of vector spaces of dimension $2$ which arises from a degeneration of linear series but is not semisimple, our Example~\ref{exac_no_SB}.

We thank Omid Amini, Marcos Jardim and Oliver Lorscheid for many discussions on the subject. 

We fix a field $k$ once and for all, and call its elements \emph{scalar}. A vector space will be a vector space over $k$. Given a vector space $V$ and $v\in V$, we let $[v]$ denote the set of all nonzero scalar multiples of $v$. We write $[v]=0$ if $v=0$. Given a map $\vp\colon V\to W$ of vector spaces we put $[\vp][v]:=[\vp(v)]$. Also, we put 
$\text{Ker}[\vp]:=\text{Ker}(\vp)$ and $\text{Im}[\vp]:=\text{Im}(\vp)$ and call them the kernel and image of $[\vp]$, respectively. We say $[\vp]$ is an isomorphism if $\vp$ is. Given a second map $\psi\colon W\to U$, we put $[\psi][\vp]:=[\psi\vp]$. 

%====================================================================================

\section{\texorpdfstring{$\mathbb Z^n$}{}-quivers}\label{Znq}

For nonnegative integers $d$ and $n$, define 
\begin{eqnarray*}
 &&\Z^{n+1}(d):= \Big\{ (d_0,\dots,d_n)\in\Z^{n+1}\,\Big|\, \sum d_i = d\Big\}.
\end{eqnarray*}
Let $v_0,\dots,v_n$ be elements in $\Z^{n+1}(0)$ such that their sum is zero and any
proper subset of them is linearly independent over $\qq$. To an element $v$ in $\Z^{n+1}(d)$ we associate a quiver $Q(v;v_0,\dots,v_n)$ whose set of vertices is
$$Q_0:=v+\Z v_0+\dots+\Z v_n \subseteq\Z^{n+1}(d)$$
and whose set of arrows $Q_1$ is the subset of $Q_0 \times Q_0$ of pairs $(u,u')$ such that
$$
u'=u+v_i \text{ for some } i\in\{0,\dots,n\}.
$$
(We view each $(u,u')$ as an arrow connecting $u$ to $u'$.) More generally, 

\begin{Def}\label{partQuiver} Let $Q$ be a quiver and $Q_0$ its set of arrows. Let $Q_0=A_0\cup\dots\cup A_n$ be a finite partition and $T:=\{A_0,\dots,A_n\}$. 
\begin{enumerate}
    \item For each $i=0,\dots,n$, call $A_i$ an \emph{arrow type} and the \emph{type} of each $a\in A_i$. 
    \item Given a path $\gamma$ in $Q$, its \emph{type} is the map $\mathfrak t_\ga\colon T\to\mathbb Z_{\geq 0}$ taking an arrow type to the number of arrows of that type $\gamma$ contains, and its \emph{essential type} is the set $\{t\in T\,|\,\mathfrak t_\ga(t)>0\}$.  
    \item A path in $Q$ is \emph{admissible} if it contains no arrows of a certain type. 
    \item A path $\ga$ in $Q$ is said to be \emph{simple} if $\mathfrak t_\ga(t)\leq 1$ for each $t\in T$.
\end{enumerate}
\end{Def}

\begin{Def}\label{defZquiver} Let $Q$ be a quiver with a nontrivial finite partition of its set of arrows into $n+1$ parts for $n>0$. We have a 
\emph{$\mathbb Z^n$-structure} on $Q$ if the following four properties are satisfied:
  \begin{enumerate}
  \item There is exactly one arrow of each type leaving each vertex.
   \item Each vertex is connected to each other by an admissible path.
  \item Each maximal simple path is a circuit.
    \item Two admissible paths leaving the same vertex arrive at the same vertex if and only if they are of the same type.
  \end{enumerate}
  A maximal simple path is called a \emph{minimal circuit} and a quiver with a $\mathbb Z^n$-structure is called a \emph{$\mathbb Z^n$-quiver}. 
  \end{Def}

Notice that the decomposition 
$Q_1=\bigcup A_i$ where 
$$
A_i:=\{(u,u') \in A \, \mid \, 
u'=u+v_i\}
$$
for $i=0,\dots,n$ gives the quiver $Q=Q(v;v_0,\dots,v_0)$ a $\mathbb{Z}^n$-structure.

\begin{Lem}\label{equivtype} 
Let $Q$ be a quiver with a nontrivial finite partition of its set of arrows in $n+1$ parts for $n>0$. Assume that there is an arrow of each type leaving each vertex. Then Properties \text{\rm (3)} and \text{\rm (4)} are equivalent to:
\begin{enumerate}
    \item[(3')] Two paths leaving the same vertex arrive at the same vertex if and only if the difference between the number of arrows of each type in each of them does not depend on the type.
\end{enumerate}
\end{Lem}

\textbf{Notation:} Given two paths $\gamma_1$ and $\gamma_2$ of a quiver $Q$ with the endpoint of $\ga_1$ equal to the initial point of $\ga_2$, their concatenation is denoted by $\ga_2\ga_1$. 
%With this, $\varphi_{\ga_2}\varphi_{\ga_1}=\varphi_{\ga_2\ga_1}$.

\begin{proof} Let $\ga_1$ and $\ga_2$ be two paths. If they are admissible or they have the same length, then they are of the same type if and only if $\mathfrak t_{\ga_1}-\mathfrak t_{\ga_2}$ is the constant function. In particular, Property (3') implies Property (4). 

Also, a maximal simple path is a path $\ga$ with $\mathfrak t_\ga$ constant, thus such that 
$\mathfrak t_\ga-\mathfrak t_\nu$ is constant, where $\nu$ is the trivial path. Then Property (3') yields that $\ga$ is a circuit.

Assume Properties (3) and (4). Then the only admissible circuit is the trivial path. Let $\gamma_1$ and $\gamma_2$ be paths leaving the same vertex. We prove they satisfy Property~(3'). Assume first they have length 2. If $\ga_1$ and $\ga_2$ are admissible, Property~(3') holds because so does Property~(4). 
Suppose, without loss of generality, that $\gamma_1$ is not admissible. Since $\gamma_1$ has length 2 and $n>0$, it follows that $n=1$ and $\gamma_1$ is a maximal simple path, whence a circuit by Property~(3). Thus $\ga_2$ arrives at the same vertex as $\ga_1$ if and only if $\ga_2$ is a circuit. Since 
$\ga_2$ has length 2 and $n=1$, it follows that $\ga_2$ is a circuit if and only if $\ga_2$ is not admissible, if and only if $\ga_2$ is a maximal simple path, if and only if $\ga_2$ has the same type as $\ga_1$, if and only if the difference between the number of arrows of each type in $\ga_1$ and $\ga_2$ does not depend on the type. So Property~(3') holds.

Now, assume $\ga_1$ and $\ga_2$ have any lengths. If $a$ and $b$ are arrows such that $ba$ is part of $\ga_1$, by hypothesis, there are an arrow $b'$ of the same type as $b$ leaving the initial vertex of $a$ and an arrow $a'$ of the same type as $a$ leaving the end vertex of $b'$. 
Then $ba$ and $a'b'$ leave from the same vertex, and arrive at the same vertex by the case we have already proved for Property~(3'). Replacing $ba$ by $a'b'$ we have a path $\ga'_1$ of the same type as $\ga_1$
% It was $\gamma$ in place of $\gamma_1$
connecting the same two vertices as $\ga_1$. 

The upshot is that we can rearrange the order of the arrows in both $\ga_1$ and $\ga_2$ and assume that $\ga_1$ and $\ga_2$ are 
sequences of maximal simple paths followed by admissible paths. As maximal simple paths are circuits by Property~(3), we may assume that $\ga_1$ and $\ga_2$ are actually admissible. But then 
Property (3') is equivalent to Property (4).
\end{proof}

\begin{proposition}\label{Znquivers}
  For each $n$, all $\mathbb Z^n$-quivers are isomorphic.
\end{proposition}

\begin{proof} Let $Q$ be the $\mathbb Z^n$-quiver $Q(v;v_0,\cdots,v_n)$. Let $Q'$ be another $\mathbb Z^n$-quiver, with vertex set $Q_0'$, arrow set $Q_1'$ and $\mathbb Z^n$-structure given by the decomposition $Q_1'=\bigcup A'_i$. It is enough to show that these two $\mathbb Z^n$-quivers are equivalent. 

We define a bijection $f_0\colon Q_0\to Q_0'$ as follows: Pick any vertex $w\in Q_0'$ and put $f_0(v):=w$. Then, for any choice of integers $\ell_0,\cdots,\ell_n$ with $\min\{\ell_i\}=0$, let
$f_0(v+\sum_i\ell_iv_i)$ be the vertex of $Q'$ obtained as the final point of a path $\gamma$ in $Q'$ with initial point $w$ and $\ell_i$ arrows of type $A'_i$ for each $i$. That $\gamma$ exists follows from Property~(1) for $Q'$, and that its final point does not depend on the chosen $\gamma$ follows from Property~(4). Now, any point of $Q_0$ is expressed as $v+\sum_i\ell_iv_i$ for unique integers $\ell_i$ with $\min\{\ell_i\}=0$. Thus $f_0$ is well-defined. It is surjective by Property~(2) for $Q'$. It is injective by Property~(4).

We define now a bijection $f_1\colon Q_1\to Q_1'$ as follows: For each $i=0,\cdots,n$ and each arrow $a\in A_i$, let $f_1(a)$ be the arrow in $A'_i$ leaving the vertex $f_0(e)$, where
$e$ is the initial point of $a$; it is unique by Property~(1) for $Q'$.  Clearly, $f_1(A_i)\subseteq A'_i$ for each $i$.
Also, $f_1$ has an inverse, mapping for each $i=0,\cdots,n$ an arrow $a'\in A'_i$ to the
arrow $a\in A_i$ connecting $e$ to $e+v_i$, where
$f_0(e)$ is the initial vertex of $a'$.

Finally, $f_1$ is compatible with $f_0$. Indeed, let $a\in Q_1$, with initial vertex $e$ and terminal vertex $h$. Then $a\in A_j$ with
$h=e+v_j$ for some $j$. Set $a':=f_1(a)$. By definition, $a'\in A_j'$ and
$u:=f_0(e)$ is its initial vertex. Let $z$ be its terminal vertex. It remains to show that $z=f_0(h)$. Now,
$$
e=v+\sum_i\ell_iv_i\quad\text{and}\quad h=v+\sum_im_iv_i
$$
for unique integers $\ell_i$ and $m_i$ with $\min\{\ell_i\}=\min\{m_i\}=0$. Also, since $u=f_0(e)$, there is a path $\gamma$ in $Q'$ connecting $w$ to $u$ having $\ell_i$ arrows of type $A'_i$ for each $i$. Then the concatenation $a'\gamma$ is a 
path in $Q'$ connecting $w$ to $z$. There are two cases:

First, assume that there is $b$ different from $j$ with $\ell_b=0$. Since $h=e+v_j$, we have that $m_i=\ell_i$ for each $i\neq j$ and $m_j=\ell_j+1$. Then $a'\gamma$ is an admissible path leaving $w$ with $m_i$ arrows of type $A'_i$ for each $i$, and hence  $z=f_0(h)$ by definition.

Finally, assume that $\ell_i>0$ for each $i$ different from $j$, whence $\ell_j=0$.
Then $m_i=\ell_i-1$ for each $i\neq j$ and $m_j=\ell_j$. It follows that $a'\gamma$ has $m_i+1$ arrows of type $A'_i$ for each $i$. By Lemma~\ref{equivtype}, any path in $Q'$ leaving $w$ with $m_i$ arrows of type $A'_i$ for each $i$ arrives at the same vertex as $a'\gamma$, which is $z$. Thus, by definition again, $z=f_0(h)$.
\end{proof}

\begin{figure}[ht]
    \centering
        \begin{subfigure}{.4\textwidth}
            \begin{tikzcd}
                \\ \\ \\ \\ \\ \\ 
                \cdots \bullet \arrow[r,shift left]
                % & \bullet \arrow[r,shift left] \arrow[l,shift left] 
                & \bullet \arrow[r,shift left]  
                \arrow[l,shift left]& \bullet \arrow[r,shift left] \arrow[l,shift left] &  \arrow[l,shift left] \bullet  \cdots
                \\
            \end{tikzcd} 
            \vspace{.2cm}
      \caption{$\Z^1$-quiver.}
    \end{subfigure}
    \hfill
    \begin{subfigure}{.4\textwidth}
        \centering
        \begin{tikzpicture}[commutative diagrams/every diagram]
            \node (p0) at (0:0cm)         {$\bullet $};
            \node (p1) at (90:1.5cm)      {$\bullet $};
            \node (p1d) at (90:2cm)       {$\rotatebox{90}{\dots}$};
            \node (p1dd) at (30:.5cm)     [shift={(90:1.5cm)}]    {$\rotatebox{30}{\dots}$};
            \node (p1ddd) at (150:.5cm)   [shift={(90:1.5cm)}]    {$\rotatebox{150}{\dots}$};
            \node (p2) at (150:1.5cm)     {$\bullet$};
            \node (p2d) at (150:.5cm)     [shift={(150:1.5cm)}] {$\rotatebox{150}{\dots}$};
            \node (p2dd) at (90:.5cm)     [shift={(150:1.5cm)}] {$\rotatebox{90}{\dots}$};
            \node (p2ddd) at (210:.5cm)   [shift={(150:1.5cm)}] {$\rotatebox{210}{\dots}$};
            \node (p3) at (210:1.5cm) {$\bullet$};
            \node (p3d) at (210:.5cm) [shift={(210:1.5cm)}] {$\rotatebox{210}{\dots}$};
            \node (p3dd) at (150:.5cm) [shift={(210:1.5cm)}] {$\rotatebox{150}{\dots}$};
            \node (p4)  at (270:1.5cm) {$\bullet $};
              \node (p5)  at (330:1.5cm) {$\bullet $};
              \node (p6)  at (30:1.5cm) {$\bullet $};
            \node (p6d) at (30:2cm) {$ \rotatebox{30}{\dots} $};
              \node (p6d) at (90:.5cm) [shift={(30:1.5cm)}]{$ \rotatebox{90}{\dots} $};
              \node (p7) at (270:3cm)   {$\bullet $};
              \node (p7d) at (270:3.5cm)   {$\rotatebox{270}{\dots}$};
            \node (p7dd) at (330:.5cm) [shift={(270:3cm)}]   {$\rotatebox{330}{\dots}$};
            \node (p7ddd) at (210:.5cm) [shift={(270:3cm)}]   {$\rotatebox{210}{\dots}$};
            \node (p8) at (0:2.6cm)           {$\bullet $};
            \node (p8d) at (30:.5cm)   [shift={(0:2.6cm)}]         {$\rotatebox{30}{\dots}$};
            \node (p8dd) at (90:.5cm)   [shift={(0:2.6cm)}]         {$\rotatebox{90}{\dots}$};
            \node (p8ddd) at (330:.5cm)   [shift={(0:2.6cm)}]         {$\rotatebox{330}{\dots}$};
            \node (p9) at (330:3cm)       {$\bullet $};
            \node (p9d) at (330:3.5cm)     {$\rotatebox{330}{\dots}$};
            \node (p9dd) at (270:.5cm)  [shift={(330:3cm)}]   {$\rotatebox{90}{\dots}$};
            \node (p9ddd) at (30:.5cm)  [shift={(330:3cm)}]   {$\rotatebox{30}{\dots}$};
            \node (p10) at (300:2.6cm)    {$\bullet$};
            \node (p10d) at (330:.5cm) [shift={(300:2.6cm)}]  {$\rotatebox{330}{\dots}$};
            \node (p10dd) at (270:.5cm) [shift={(300:2.6cm)}]  {$\rotatebox{90}{\dots}$};
            \node (p11) at (240:2.6cm)    {$\bullet$};
            \node (p11d) at (210:.5cm)  [shift={(240:2.6cm)}] {$\rotatebox{210}{\dots}$};
            \node (p11dd) at (270:.5cm)  [shift={(240:2.6cm)}] {$\rotatebox{90}{\dots}$};
            \node (p11ddd) at (150:.5cm)  [shift={(240:2.6cm)}] {$\rotatebox{150}{\dots}$};
            \path[commutative diagrams/.cd, every arrow, every label]  
            (p8) edge[<-] (p5)    (p9) edge[<-]  (p8)
            (p10) edge[<-] (p5)   (p9) edge[<-] (p10)
            (p4) edge[<-] (p10)   (p10) edge[<-] (p7)
            (p0) edge[<-] (p1)    (p2) edge[<-] (p0)
            (p0) edge[<-] (p3)    (p4) edge[<-] (p0) 
            (p0) edge[<-] (p5)    (p6) edge[<-] (p0) 
            (p6) edge[<-] (p8)    (p5) edge[<-] (p9)
            (p7) edge[<-] (p4)    (p11) edge[<-] (p3) 
            (p4) edge[<-] (p11)   (p11) edge[<-] (p7) 
  
            (p1) edge[<-] (p2)    (p1) edge[<-] (p6) 
            (p5) edge[<-] (p6)    (p5) edge[<-] (p4) 
            (p3) edge[<-] (p4)    (p3) edge[<-] (p2) ;
        \end{tikzpicture}
        \vspace{-.3cm}
        \caption{$\Z^2$-quiver.}
        \end{subfigure}
    \caption{$\mathbb Z^n$-quivers}
    \label{fig1}
\end{figure}
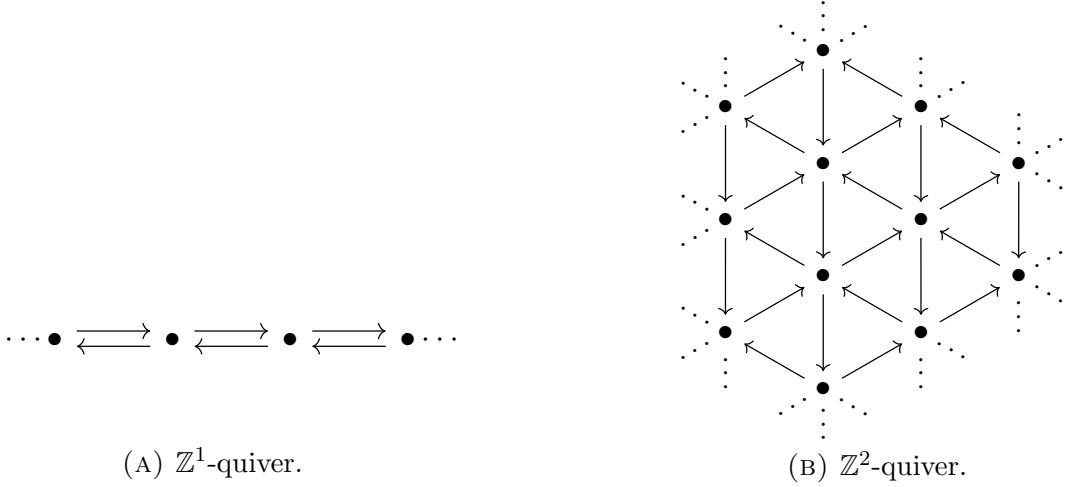

We draw in Figure-\ref{fig1} the $\mathbb{Z}^1$-quiver and the $\mathbb{Z}^2$-quiver, up to isomorphism.

\section{Degenerations of linear series}
\label{lls_qui_rep_sb}

$\mathbb{Z}^n$-quivers arise from degenerations of linear series. 
Indeed, let $X$ be a reduced and reducible scheme of finite type over the fixed field $k$.  
%Assume its irreducible components are geometrically irreducible; 
Order and denote the irreducible components by $X_0,\dots,X_n$ for $n> 0$. Let $B$ be the spectrum of a discrete valuation ring $R$ with residue field $k$. Let 
$\pi\colon\mathcal{X}\to B$ be a flat morphism with special fiber isomorphic to $X$. 
%We call it a \emph{smoothing} of $X$ if its generic fiber is smooth. 
We say $\pi$ is
\emph{regular} if $\mathcal{X}$ is. We identify the special fiber with $X$. In what follows we assume that $\pi$ is regular. 
%(A nonregular smoothing of a projective curve can be made into a regular one by the process known as semistable reduction, at the cost of replacing the special fiber.)

Let $L_\eta$ be an invertible sheaf on the generic fiber $\mathcal{X}_\eta$ of $\pi$. 
As $\mathcal{X}$ is regular, each $X_i$ is a Cartier divisor of $\mathcal X$ and there exists an invertible sheaf extension $\mathcal{L}$ of $L_\eta$ to $\mathcal{X}$. The extension is not unique: For each divisor $D=\sum\ell_iX_i$, the invertible sheaf $\mathcal{L}(D):=\mathcal{L}\otimes\mo_\mathcal{X}(D)$ is also an extension of
$L_\eta$. 

Fix an invertible extension $\mathcal L$ of $L_\eta$. The $\mathcal L(D)$ are all the invertible extensions of $L_\eta$. Indeed, given an invertible extension $\mathcal M$, there is a map $\mathcal L\to\mathcal M$ extending the identity map on $L_\eta$. Thus, there is an effective Cartier divisor $D$ on $\mathcal X$ supported on $X$ such that $\mathcal M\cong\mathcal L(D)$. The corresponding Weil divisor is of the form $\sum\ell_iX_i$. Since $\mathcal X$ is regular, the latter is also a Cartier divisor and coincides with $D$.

To each $u=(\ell_0,\dots,\ell_n)\in\mathbb Z^{n+1}$ associate the invertible sheaf $\mathcal L_u:=\mathcal L(\sum\ell_iX_i)$ on $\mathcal X$ and its restriction $L_u:=\mathcal L_u|_X$ to $X$. Let $u_{-1}:=u-(1,\dots,1)$. Since $X$ is reduced, the natural map $\mathcal L_{u_{-1}}\to\mathcal L_u$ factors through an isomorphism $\mathcal L_{u_{-1}}\to T\mathcal L_u$, allowing us to identify $\mathcal L_{u_{-1}}$ with $\mathcal L_u$.

Let $Q_0\subseteq\mathbb Z^{n+1}$ be the set of $(n+1)$-tuples $(\ell_0,\dots,\ell_n)$ of nonnegative integers such that $\min\{\ell_i\}=0$. 
For each $u\in\mathbb Z^{n+1}$ put $\overline u:=u-\min\{u_i\}(1,\dots,1)$. Then $\overline u\in Q_0$ and $\mathcal L_u\cong\mathcal L_{\overline u}$.

Given $u_1,u_2\in Q_0$, let $D^{u_1}_{u_2}=\sum\ell_iX_i$, where 
$(\ell_0,\dots,\ell_n):=\overline{u_2-u_1}$. Twisting by $D^{u_1}_{u_2}$ and identifying $\mathcal L_{u_2}$ with $\mathcal L_{u_1+\overline{u_2-u_1}}$ yields a map $\psi^{u_1}_{u_2}\colon\mathcal L_{u_1}\to\mathcal L_{u_2}$ which restricts to an isomorphism on $\mathcal X_{\eta}$ and to a nonzero map $\lambda^{u_1}_{u_2}\colon L_{u_1}\to L_{u_2}$ on $X$. In fact, 
\begin{equation}\label{prop3}
\lambda^{u_1}_{u_2}|_{X_i}=0\quad \text{if and only if}\quad \ell_i>0
\end{equation}
for each $i=0,\dots,n$. 

Given $u_1,u_2,u_3\in Q_0$ we have that $D^{u_2}_{u_3}+D^{u_1}_{u_2}$ does not contain the whole $X$ in its support if and only if $\overline{u_3-u_2}+\overline{u_2-u_1}\in Q_0$, in which case $\psi_{u_3}^{u_1} = \psi_{u_3}^{u_2} \circ \psi_{u_2}^{u_1}$. It follows that 
\begin{equation}\label{comp}
\lambda_{u_3}^{u_2} \circ \lambda_{u_2}^{u_1}=
\begin{cases}
\lambda_{u_3}^{u_1}&\text{if } \overline{u_3-u_2}+\overline{u_2-u_1}\in Q_0,\\
0&\text{otherwise.}
\end{cases}
\end{equation}

For each $i=0,\dots,n$, let $\chi_i\in\mathbb Z^{n+1}$ be the tuple whose $j$-th component is $0$ for $j\neq i$ and whose $i$-th component is $1$, and let $A_i\subseteq Q_0\times Q_0$ be the subset of pairs $(u_1,u_2)$ such that $u_2=\overline{u_1+\chi_i}$. Then the $A_i$ are pairwise disjoint. Put $Q_1:=\bigcup A_i$. We obtain a quiver $Q:=(Q_0,Q_1)$. Furthermore, the $A_i$ form a nontrivial partition of $Q_1$ giving $Q$ a $\mathbb Z^n$-structure. 

In addition, the $L_u$ for $u\in Q_0$ and the maps $\lambda^{u_1}_{u_2}$ for 
$(u_1,u_2)\in Q_1$ give us a representation $\mathfrak L$ of $Q$ in the category of invertible sheaves on $X$. It is a special representation.

\begin{proposition}\label{Lgen} $\mathfrak L$ is a \emph{maximal exact linked net} over $Q$ of invertible sheaves on $X$.
\end{proposition}

\begin{proof} The definitions and proof will be given in Section~\ref{lnets}.
\end{proof}

Now, let $V_{\eta}\subseteq H^0(\mathcal{X}_\eta,\mathcal{L}_\eta)$ be a vector subspace. Thus $(L_\eta,V_\eta)$ is a linear series on $\mathcal X_\eta$. For each $u\in Q_0$ 
there is an extension $\mv_{u}$ of $V_\eta$ given by 
$$
\mv_u:=\Big\{s\in H^0(\mathfrak{X},\mathcal{L}_{u})\,\Big|\, s|_{\mathfrak{X}_\eta}\in \mv_\eta\Big\}.
$$
Set $V_{u}:={\mv_{u}}|_X\subseteq H^0(X,L_{u})$, which is a vector space over $k$. For each $u_1,u_2\in Q_0$ the map $\psi^{u_1}_{u_2}$ induces maps (abusing notation)
$$
\psi^{u_1}_{u_2} \colon \mv_{u_1} \longrightarrow \mv_{u_2}\quad\text{and}\quad \lambda^{u_1}_{u_2} \colon V_{u_1} \to V_{u_2}.
$$

The $\mathcal V_u$ for $u\in Q_0$ and the maps $\psi^{u_1}_{u_2}$ for $(u_1,u_2)\in Q_1$ give us a representation $\mathfrak M$ of $Q$ in the category of modules over $R$. Its restriction to $X$ (resp.~$\mathcal X_{\eta}$) is a representation of $Q$ in the category of vector spaces over $k$ (resp.~the field of fractions $K$ of $R$). The latter has associated maps which are isomorphisms, 
whereas the former is a subrepresentation $\g$ of $H^0(X,\mathfrak L)$, the induced representation of vector spaces obtained by taking global sections in $\mathfrak L$.

Assume $\pi$ is proper. Then the spaces $V_u$ are finite-dimensional. Their dimension is the same as the rank of the $\mathcal V_u$, which is the dimension of $V_\eta$ over the field $K$. Let $H$ be the collection of $u\in Q_0$ for which there is a section of $L_u$ in $V_u$ that does not vanish on any irreducible component of $X$. We call $H$ the \emph{effective locus} of $(\mathfrak L,\g)$.

\begin{proposition}\label{Vgen} $H^0(X,\mathfrak L)$ and $\g$ are \emph{linked nets} over $Q$ of vector spaces over $k$, and $\g$ is \emph{exact}. If $\pi$ is proper (and thus $X$ is complete over $k$) then $\g$ is \emph{pure}. If in addition $X$ is connected, then $H$ is finite and $\g$ is \emph{$1$-generated} by $H$.
\end{proposition}

\begin{proof} The definitions and proof will be given in Section~\ref{lnets}.
\end{proof}

The quiver we are considering here may look different from that in \cite{Osserman_2019}. It is only apparently so, as we now explain.

Assume now that $X$ is a \emph{curve}, that is, $X$ has irreducible components of dimension 1. Assume as well that $X$ is connected. 

For each $i,j=0,\dots,n$, let $x_i:=\deg\mathcal{L}|_{X_i}$ and 
$x_{i,j}:=\deg\mo_{\mathcal{X}}(X_i)|_{X_j}$. Put $v=(x_0,\dots,x_n)$ and $v_i:=(x_{i,0},\dots,x_{i,n})$ for $i=0,\dots,n$. Then $v\in\mathbb Z^{n+1}(d)$ and 
$v_i\in\mathbb Z^{n+1}(0)$ for each $i$, where $d$ is the degree of $L_\eta$. 
Since $X$ is reduced, $\mathcal O_{\mathcal X}(\sum X_i)\cong\mathcal O_{\mathcal X}$, and thus $\sum v_i=0$. Also, since $X$ is connected, $x_{i,j}< 0$ if and only if $i=j$, and thus each proper subset of the $v_i$ is linearly independent 
over $\mathbb Q$.

Put 
$$
Q_0':=v+\mathbb Zv_0+\cdots+\mathbb Zv_n.
$$
The natural function $f_0\colon Q_0\to Q_0'$ defined by letting 
$f_0(\ell_0,\dots,\ell_n):=v+\sum\ell_iv_i$ is thus a bijection. Furthermore, 
for each $u=(\ell_0,\dots,\ell_n)\in\mathbb Z^{n+1}$, the multidegree of $L_u$ is $v+\sum\ell_iv_i$. It follows that the multidegrees of $L_{u_1}$ and $L_{u_2}$ are distinct for distinct $u_1,u_2\in Q_0$.

%Since $X$ is connected, $L_u\cong\mathcal L_v$ if and only if $\overline{u-v}=0$. Indeed, it is enough to observe that $\mathcal O_{\mathfrak X}(\sum\ell_i X_i)|_X\cong\mathcal O_{X}$ if and only if $(\ell_0,\dots,\ell_n)\in\mathbb Z(1,\dots, 1)$. Since $\mathcal O_{\mathfrak X}(\sum X_i)\cong\mathcal O_{\mathfrak X}$ we may assume $\min\{\ell_i\}=0$. Without loss of generality, assume 

For each $i=0,\dots,n$, let $ A'_i \subset Q_0' \times Q_0'$ be the subset 
of pairs $(w_1,w_2)$ such that $w_2=w_1+v_i$,  and set $Q_1':=\bigcup A'_i$. The bijection $f_0$ induces a bijection 
$f_1\colon Q_0\times Q_0\to Q_0'\times Q_0'$ taking $A_i$ onto $A'_i$ for $i=0,\dots,n$. Thus 
$f_0$ induces an isomorphism between the quiver $Q$ and the quiver
$Q(v;v_0,\dots,v_n)$. The latter is the quiver that appears in \cite{Osserman_2019}, Notation~2.11, p.~62, for the trivial chain structure; see loc.~cit., Def.~2.2, p.~60.

Notice that $f_0$ takes the effective locus $H$ of $(\mathfrak L,\g)$ to a subset of the set of effective multidegrees, that is,
$\mathbb Z^{n+1}_{\geq 0}\cap\mathbb Z^{n+1}(d)$, which is clearly finite.

%=================================================================================

\section{Linked nets}\label{lnets}

Given a representation $\g$ of a quiver $Q$ in an $k$-linear Abelian category, the object associated to a vertex $v$ of $Q$ will be denoted $V^\g_v$, and the composition of maps associated to a path $\gamma$ in $Q$ will be denoted $\vp^\g_\gamma$. (If $\gamma$ is the trivial path, $\vp^\g_\gamma=\text{id}$ by convention.) We drop the superscript $\g$ if $\g$ is clear from the context. 

%From now on, unless stated otherwise, a representation $\g$ of a quiver $Q$ is one in the category of finite dimensional vector spaces over a fixed field. We denote by $V^{\mathfrak g}_v$ the vector space associated to each vertex $v$ of $Q$ and by $\varphi^{\mathfrak g}_{\gamma}$ the composition of maps associated to each path $\gamma$ in $Q$. We drop the superscript $\g$ if $\g$ is clear from the context. An element of $V^{\mathfrak g}_v$ will be called a \emph{section of $\mathfrak g$ at $v$.}

\begin{Def}\label{deflnvsfs} Let $\g$ be a representation of a quiver $Q$ in a $k$-linear Abelian category. Let $H$ be a set of vertices of $Q$ and $m\in\mathbb N\cup\{\infty\}$. We say:
\begin{enumerate} 
\item $\g$ is \emph{pure} if every epimorphism between the objects of the category associated to the vertices is an isomorphism.
\item $\g$ is \emph{locally finite} if for each vertex $v$ of $Q$ there is an integer $\ell$ such that $\vp_\mu=0$ for every path $\mu$ arriving at $v$ of length greater than $\ell$.
\item $\g$ is (\emph{faithfully}) \emph{$m$-generated} by $H$ if for each vertex $v$ of $Q$ there exist $w_1,\dots,w_n\in H$ for $n\in\mathbb N$ with $n\leq m$ and paths $\mu_1,\dots,\mu_n$ connecting them to $v$ such that 
$$
(\vp_{\mu_1},\cdots,\vp_{\mu_n})\colon\bigoplus_{i=1}^m V_{w_i}\longrightarrow V_v
$$
is an epimorphism (an isomorphism).
\end{enumerate}
If $\g$ is $\infty$-generated by $H$, we simply say $\g$ is generated by $H$. And if in addition $H$ is finite, we say $\g$ is \emph{finitely generated}.
\end{Def}

If $\g$ is pure and $1$-generated by $H$ then $\g$ is faithfully $1$-generated by $H$.

\begin{Def}\label{deflnvs} Let $\g$ be a representation of a $\mathbb Z^n$-quiver $Q$ in a $k$-linear Abelian category. 
We call $\g$ a 
\emph{weakly linked net} over $Q$ (of objects in the category) if:
  \begin{enumerate}
  \item If $\gamma_1$ and $\gamma_2$ are two paths connecting the same two vertices and $\gamma_2$ is admissible then $\varphi_{\gamma_1}$ is a scalar multiple of $\vp_{\gamma_2}$. 
  \item $\vp_\gamma=0$ for each minimal circuit $\gamma$ of $Q$. 
 \end{enumerate}
We call $\g$ a \emph{linked net} if furthermore:
  \begin{enumerate}
  \item[(3)] If $\gamma_1$ and $\gamma_2$ are two admissible paths leaving the same vertex with no arrow type in common then
    $\text{Ker}(\varphi_{\gamma_1})\cap \text{Ker}(\varphi_{\gamma_2})=0$.
    \end{enumerate}
\end{Def}

Let $\g$ be a weakly linked net over a $\mathbb Z^n$-quiver $Q$. For each two vertices $u_1$ and $u_2$ of $Q$, put $\varphi_{u_2}^{u_1}:=[\varphi_\gamma]$ for an admissible path $\gamma$ connecting $u_1$ to $u_2$. Property~(1) for $\g$ yields that $\varphi_{u_2}^{u_1}$ is well-defined. 

%As above, none of the statements we will make involving the $\varphi_{u_2}^{u_1}$ depend on the choice of $\gamma$. 

\begin{Def}\label{neighbors} Distinct vertices $u_1$ and $u_2$ of $Q$ are called \emph{neighbors} if there is an admissible simple path $\gamma$ connecting $u_1$ to $u_2$. 
\end{Def}

If $u_1$ and $u_2$ are neighbors, there is also an admissible simple path $\mu$ connecting $u_2$ to $u_1$, and both concatenations $\mu\gamma$ and $\gamma\mu$ are minimal circuits. It follows from Property~(2) that $\varphi^{u_2}_{u_1}\varphi^{u_1}_{u_2}=0$.

\begin{Def}\label{exact}
  A weakly linked net $\g$ over a $\mathbb{Z}^n$-quiver $Q$ is \emph{exact} if
  $$
  \text{Ker}(\varphi^{u_1}_{u_2})=\text{Im}(\varphi^{u_2}_{u_1})
  $$
  for each two neighboring vertices $u_1$ and $u_2$ of $Q$.
\end{Def}

%\begin{Def}\label{deflnvsfs} Let $\g$ be a weakly linked net of vector spaces, and $H$ a set of vertices of $Q$. We say $\g$ is \emph{supported} on $H$ if for each vertex $w$ of $Q$ there exists a path $\gamma$ connecting $H$ to $w$ such that $\varphi_{\gamma}$ is an isomorphism. If in addition $H$ is finite we say $\mathfrak g$ has \emph{finite support}. We say $\g$ is \emph{generated} by $H$ if for each vertex $v$ of $Q$, the space $V_v$ is the sum of $\text{Im}(\varphi_\gamma)$ for all paths $\gamma$ connecting $H$ to $v$. If in addition $H$ is finite, we say $\mathfrak g$ is \emph{finitely generated}.
%\end{Def}

%A general linked net has support on and is generated by each vertex of the quiver. 

%If $H$ is a subset of vertices as in the previous definition, we say that $\mathfrak{g}$ has \emph{finite support on} $H$. 

%\begin{Rem}
% A linked net arising from a degeneration of a linear series has finite support.
%\end{Rem}

%An exact weakly linked net of vector spaces is special: each minimal circuit $\gamma$ is a concatenation $\gamma=\gamma_1\gamma_2$ of two nontrivial simple admissible paths $\gamma_1,\gamma_2$, and thus $\varphi_{\gamma_2}\varphi_{\gamma_1}=0$ by exactness.

\begin{Def}\label{maximal} Let $X$ be a reduced scheme of finite type over $k$. A weakly linked net $\mathfrak L$ over a $\mathbb{Z}^n$-quiver $Q$ in the category of invertible sheaves on $X$ is said to be \emph{maximal} if the map associated to each arrow of $Q$ is generically zero on one and only one irreducible component of $X$.
\end{Def}

%\begin{proposition} Let $\mathfrak L$ be the representation in the category of line bundles over $X$, and $\g$ the subrepresentation of the representation by vector spaces $H^0(X,\mathfrak L)$ obtained from $\mathfrak L$ by taking global sections, of the quiver arising from a degeneration of linear series. Then $\mathfrak L$, $H^0(\mathfrak L)$ and $\g$ are linked nets, the first exact and maximal, the third pure, finitely generated and exact.
%\end{proposition}

{\noindent {\bf Proof of Proposition~\ref{Lgen}.} Recall the notation used in Section~\ref{lls_qui_rep_sb}. That $\mathfrak L$ is a linked net follows from the definition: In fact, the two properties of a weakly linked net follow from \eqref{comp}, whereas the third property of a linked net follows from \eqref{prop3} and the fact that $X$ is reduced.}

Furthermore, \eqref{prop3} implies $\mathfrak L$ is maximal. It is also exact. Indeed, let $u_1$ and $u_2$ be neighboring vertices of $Q$. Then 
\begin{equation}\label{Lu1u2}
\mathcal L_{u_2}=\mathcal L_{u_1}\big(\sum_{i\in I}X_i\big)\quad\text{and}\quad
\mathcal L_{u_1}=\mathcal L_{u_2}\big(\sum_{i\not\in I}X_i\big)
\end{equation}
for a certain proper nonempty subset $I\subset\{0,\dots,n\}$.
Set $Y:=\bigcup_{i\in I}X_i$ and $Z:=\bigcup_{i\not\in I}X_i$ with the reduced induced structures. Then the maps $\psi^{u_1}_{u_2}$ and $\psi^{u_2}_{u_1}$ induce right-exact sequences upon restriction to $X$:
$$
\begin{CD}
L_{u_1} @>\lambda^{u_1}_{u_2}>> L_{u_2} @>>> L_{u_2}|_Y @>>> 0,\\
L_{u_2} @>\lambda^{u_2}_{u_1}>> L_{u_1} @>>> L_{u_1}|_Z @>>> 0.
\end{CD}
$$
Because of \eqref{prop3}, since $X$ is reduced, 
$\lambda^{u_1}_{u_2}$ factors through an injection $L_{u_1}|_Z\hookrightarrow L_{u_2}$, and $\lambda^{u_2}_{u_1}$ factors through an injection $L_{u_2}|_Y\hookrightarrow L_{u_1}$. The exactness of $\mathfrak L$ follows.\qed

\

{\noindent {\bf Proof of Proposition~\ref{Vgen}.} Recall the notation used in Section~\ref{lls_qui_rep_sb}. That $\mathfrak L$ is a linked net follows from Proposition~\ref{Lgen}. But then it follows that also $H^0(X,\mathfrak L)$ and $\g$ are linked nets.}

Also, $\mathfrak V$ is exact. Indeed, let $u_1$ and $u_2$ be neighboring vertices, and $s\in V_{u_1}$ such that $\lambda^{u_1}_{u_2}(s)=0$. Then $s=\widetilde s|_X$ for some $\widetilde s\in\mathcal V_{u_1}$. Since $u_1$ and $u_2$ are neighbors, there is a nonempty proper subset $I\subset\{0,\dots,n\}$ such that \eqref{Lu1u2} holds.  For each $i\not\in I$, since $\psi^{u_1}_{u_2}$ is generically an isomorphism on $X_i$, we have $\widetilde s|_{X_i}=0$. Hence $\widetilde s=\psi^{u_2}_{u_1}(\widetilde x)$ for some $\widetilde x\in H^0(\mathcal X,\mathcal L_{u_2})$. Since $\psi^{u_2}_{u_1}$ restricts to the identity on the generic fiber, $\widetilde s$ and $\widetilde x$ agree there, and thus $\widetilde x\in\mathcal V_{u_2}$. It follows that $x:=\widetilde x|_X$ is in $V_{u_2}$ and $s=\lambda^{u_2}_{u_1}(x)$.

Assume $\pi$ is a proper map. Then $\mathcal V_u$ is a free $R$-module of finite rank equal to the dimension of $V_\eta$ over $K$, and thus $\g$ has pure finite dimension; in fact, all the spaces $V_u$ have dimension over $k$ equal to the dimension of $V_\eta$ over $K$. Thus $\g$ is pure.

In addition, assume $X$ is connected. Then the effective locus $H$ of $(\mathfrak L,\mathfrak V)$ is finite. Indeed, let $u_\ell=(u_{\ell,0},\dots,u_{\ell,n})\in Q_0$ forming an infinite sequence of distinct elements. We need only prove that there is $\ell$ such that $u_\ell\not\in H$. Since for each $\ell$ there is $i\in\{0,\dots,n\}$ such that $u_{\ell,i}=0$, up to replacing the sequence $u_\ell$ by an infinite subsequence, we may assume that there is $i$ such that $u_{\ell,i}=0$ for each $\ell$. Let 
$$
I:=\Big\{i\in\{0,\dots,n\}\,|\,u_{\ell,i} \text{ is unbounded}\Big\}.
$$
Then $I$ is a nonempty proper subset of $\{0,\dots,n\}$. Up to replacing the sequence $u_\ell$ by an infinite subsequence, we may also assume that $u_{\ell,i}\geq\ell+m$ for each $\ell$ and $i\in I$ and $u_{\ell,j}\leq m$ for each $\ell$ and $j\not\in I$, for a certain $m\in\mathbb N$. Since $X$ is connected, there is $i\in I$ and $j\not\in I$ such that $X_i\cap X_j\neq\emptyset$. Since $X_i$ and $X_j$ are Cartier divisors of $\mathcal X$, the intersection $X_i\cap X_j$ is a Cartier divisor of $X_i$. Then 
$$
\mathcal L_{u_\ell}|_{X_i}\subseteq\mathcal L_0|_{X_i}\Big(-\ell X_i\cap X_j\Big)
$$
for each $\ell$. Since $X_i$ is integral, no nonzero section of $\mathcal L_0|_{X_i}$ can vanish to infinite order on $X_i\cap X_j$ for any $\ell$. Furthermore, since 
$\pi$ is proper, $X_i$ is complete and thus $H^0(X_i,\mathcal L_0|_{X_i})$ is finite dimensional. It follows that there is $\ell$ such that 
$H^0(X_i,\mathcal L_{u_\ell}|_{X_i})=0$. But then, every section of $\mathcal L_{u_\ell}$ vanishes on $X_i$ and thus $u_\ell\not\in H$. 

Observe that, since $\pi$ is regular, $\mathcal X$ is reduced. It is also irreducible if connected. But, since $X$ is connected, so is $\mathcal X$, because $\pi$ is proper and the special point of $B$ is its unique nonempty proper closed subset. So $\mathcal X$ is integral. We claim $\g$ is $1$-generated by $H$. Indeed, let $u\in Q_0-H$. Then each section of $L_u$ in $V_u$ vanishes on $X_i$ for some $i$. It follows that there is $i$ such that every section of $L_u$ in $V_u$ vanishes on $X_i$, and hence every section of $\mathcal L_u$ in $\mathcal V_u$ vanishes on $X_i$. Then $\psi^{z}_u\colon\mathcal V_z\to\mathcal V_u$ is surjective, where $z:=\overline{u-\chi_j}$. If $z\not\in H$, repeat the process. Since $\mathcal X$ is integral, it must end because a nonzero global section of $\mathcal L_u$ vanishes to finite order on each component $X_i$. We will have thus found a sequence $z_0,\dots,z_m\in Q_0$ with $z_0=u$ and $z_m\in H$ such that $\psi^{z_i}_{z_{i-1}}\colon\mathcal V_{z_i}\to\mathcal V_{z_{i-1}}$ is surjective and hence so is $\lambda^{z_i}_{z_{i-1}}\colon V_{z_i}\to V_{z_{i-1}}$ for each $i=1,\dots,m$. But then it follows from \eqref{comp} that $\lambda^{z_m}_u$ is surjective. 
\qed

\section{Properties of \texorpdfstring{$\Z^n$}{}-quivers}\label{generalquiver}

Let $Q$ be a $\mathbb Z^n$-quiver. Let $\ga$ be an admissible path in $Q$. Then there is an admissible  path $\overline\ga$ in $Q$ such that $\ga\overline\ga$ and $\overline\ga\ga$ make sense and have the same type as concatenations of minimal circuits. We call $\overline\ga$ a \emph{reverse path} to $\ga$. It connects the end vertex of $\ga$ to its inital vertex.

\begin{proposition} Let $Q$ be a $\mathbb Z^n$-quiver. Let $\ga$ be an admissible path in $Q$ and $\overline\ga$ a reverse path. Then the following statements hold:
\begin{enumerate}
    \item $\ga$ is trivial if and only if $\overline\ga$ is trivial.
  %  \item If $\ga$ is nontrivial, a type of arrow appears in $\ga$ if and only if it does not appear in $\ol\ga$.
    \item $\ga$ is simple if and only if $\overline\ga$ is simple.
    \item $\ga$ is a reverse path to $\overline\ga$. 
\end{enumerate}
\end{proposition}

\begin{proof} All of the statements are straightforward.
\end{proof}

\begin{Lem}\label{adjadm} Let $Q$ be a $\mathbb{Z}^n$-quiver with vertex set $Q_0$. Let $u,v,w\in Q_0$. Then the following first two statements are equivalent and imply the last two:
  \begin{enumerate}
  \item There is an admissible path from $u$ to $w$ through $v$.
  \item The concatenation of every admissible path from $u$ to $v$ and every admissible path from $v$ to $w$ is admissible.
  \item If $v\neq w$ then all paths from $u$ to $v$ through $w$ are not admissible.
  \item If $u\neq v$ then all paths from $v$ to $w$ through $u$ are not admissible.
    \end{enumerate}
  Furthermore, if $v$ and $w$ are neighbors then the third statement is equivalent to the first two, and if $u$ and $v$ are neighbors then the fourth statement is equivalent to the first two.
\end{Lem}

\begin{proof} Assume Statement (1). Let $\gamma$ be an admissible path from $u$ to $w$ through $v$. Then $\ga=\ga_2\ga_1$ for a path $\gamma_1$ from $u$ to $v$ and a path $\gamma_2$ from $v$ to $w$, both admissible. If $\beta$ is an admissible path from $u$ to $v$ and $\nu$ is one from $v$ to $w$, then $\beta$ has the same type as $\ga_1$ and $\nu$ has the same type as $\ga_2$, whence $\nu\beta$ has the same type as $\ga_2\ga_1$. 
Since $\ga=\ga_2\ga_1$, the concatenation $\nu\beta$ is admissible because so is $\ga$. 

Assume Statement~(2) now. Since there are admissible paths from $u$ to $v$ and from $v$ to $w$, Statement (1) holds. Furthermore, assume $v\neq w$. Let $\gamma$ be a path from $u$ to $w$ and $\mu$ one from $w$ to $v$. Assume by contradiction to Statement (3) that $\beta:=\mu\ga$ is admissible. Then $\gamma$ and $\mu$ are admissible. Also, a type of arrow not appearing in $\beta$ does not appear in $\mu$ either. Let $\nu$ be a reverse path to $\mu$. Since $v\neq w$, a type of arrow not appearing in $\mu$ appears in $\nu$. It follows that $\nu\beta$ is not admissible, contradicting Statement (2).

Assume Statement (2) again. Assume $u\neq v$. Let $\gamma$ be a path from $v$ to $u$ and $\mu$ one from $u$ to $w$. Assume by contradiction to Statement (4) that $\beta:=\mu\ga$ is admissible. Then $\gamma$ and $\mu$ are admissible. Also, a type of arrow not appearing in $\beta$ does not appear in $\ga$ either. Let $\nu$ be a reverse path to $\gamma$. Since $u\neq v$, a type of arrow not appearing in $\ga$ appears in $\nu$.  It follows that $\beta\nu$ is not admissible, negating Statement (2).
% \textcolor{red}{This argument is very similar to the previous one, maybe it can be omitted?}

Let now $\beta$ be an admissible path from $u$ to $v$ and $\nu$ one from $v$ to $w$. Assume $\nu\beta$ is not admissible, negating Statement~(2). Let $\ga$ be an admissible path connecting $u$ to $w$ and $\epsilon$ a minimal circuit connecting $w$ to itself. If $\beta$ or $\nu$ is simple, then $\epsilon\ga$ has the same type as $\nu\beta$.
  
If $v$ and $w$ are neighbors (thus $v\neq w$), then $\nu$ is simple. So is a reverse path $\mu$. Furthermore, since 
$\mu\nu\beta$ and $\mu\epsilon\ga$ have the same type, so have $\beta$ and $\mu\ga$. Since $\beta$ is admissible, so is $\mu\ga$, negating Statement~(3).

If $u$ and $v$ are neighbors (thus $u\neq v$), then $\beta$ is simple. So is a reverse path $\alpha$. Furthermore, since 
$\nu\beta\alpha$ and $\epsilon\ga\alpha$ have the same type, so have $\nu$ and $\ga\alpha$. Since $\nu$ is admissible, so is $\ga\alpha$, negating Statement~(4).
% \textcolor{red}{This argument is very similar to the previous one, maybe it can be omitted?}
\end{proof}  

\begin{Def} Let $Q$ be a $\Z^n$-quiver. For each vertex $v$ of $Q$ and each set $I$ of arrow types, the \emph{$I$-cone of $v$}, denoted $C_I(v)$, is the set of end vertices of all the paths leaving $v$ having essential types contained in $I$.
\end{Def}

\begin{proposition}\label{Cfinite}
Let $Q$ be a $\mathbb Z^n$-quiver, and $z_0,\dots,z_m$ be vertices of $Q$. Let $I_0,\dots,I_m$ be sets of arrow types such that $I_0\cap\cdots\cap I_m=\emptyset$. Then 
$$
   \bigcap_{j=0}^m C_{I_j}(z_j) 
$$
is finite.
\end{proposition}

\begin{proof} Setting $I_i:=I_m$ and $z_i:=z_m$ for $i=m+1,\cdots,n$ we may assume $m=n$ and $\# I_j=n$ for each $j$. Since all $\mathbb Z^n$-quivers are equivalent by Proposition~\ref{Znquivers}, we may assume that $Q$ is
the quiver $Q(v;v_0,\dots,v_n)$ with $v=0$, with $v_1,\dots,v_n$ forming the canonical basis of $\mathbb Z^n$ and $v_0=(-1,\dots,-1)$. We may also assume that $A_j$ is the only type not in $I_j$ for each $j$. Thus
\begin{align*}
C_{I_j}(z_j)=&\Big\{x\in\mathbb Z^n\,\Big|\,x_j\leq z_{j,j}\Big\}\quad\text{for }j=1,\dots,n,\\
C_{I_0}(z_0)=&\Big\{x\in\mathbb Z^n\,\Big|\,x_j\geq z_{0,j}\text{ for }j=1,\dots,n\Big\}.
\end{align*}
Their intersection is thus the finite set
$$
[z_{0,1},z_{1,1}]\times [z_{0,2},z_{2,2}]\times\cdots\times
[z_{0,n},z_{n,n}].
$$
\end{proof}

\begin{Def}\label{hull} Let $Q$ be a $\mathbb Z^n$-quiver. Let $H$ be a non-empty set of vertices of $Q$. Let $P(H)$ be the set of all vertices $v$ of $Q$ such that for each arrow type there are $z\in H$ and a path $\gamma$ connecting $z$ to $v$ and containing no arrow of that type. We call $P(H)$ the \emph{hull} of $H$.
\end{Def}

\begin{proposition}\label{HPH} Let $Q$ be a $\mathbb Z^n$-quiver.
  Let $H$ be a non-empty set of vertices of $Q$. Then the following three statements hold:
  \begin{enumerate}
  \item $H\subseteq P(H)$.
  \item If $H$ is finite, so is $P(H)$.
    \item $P(P(H))=P(H)$.
  \end{enumerate}
\end{proposition}

\begin{proof}
Statement (1) is clear: If $v\in H$ then the trivial path connects a vertex of $H$ to $v$ and contains no arrows of any type.

As for Statement~(2), observe that
$$
  P(H)=\bigcup_{f\in H^T} \bigcap_{\mathfrak a\in T} C_{I_{\mathfrak a}}(f(\mathfrak a)),
$$
where $T$ is the set of arrow types and $I_{\mathfrak a}:=T-\{\mathfrak a\}$ for each $\mathfrak a\in T$.
If $H$ is finite, so is the set of functions $H^T$. Thus $P(H)$ is a finite union of
finite sets by Proposition~\ref{Cfinite}.

As for Statement~(3), the inclusion $P(H)\subseteq P(P(H))$ follows from Statement~(1). In addition, for each $v\in P(P(H))$ and $\mathfrak a\in T$
there are $w_{\mathfrak a}\in P(H)$ and a path $\gamma_{\mathfrak a}$ connecting $w_{\mathfrak a}$ to $v$ and containing no arrow of type $\mathfrak a$. Since $w_{\mathfrak a}\in P(H)$, there are $z_{\mathfrak a}\in H$ and a path
$\mu_{\mathfrak a}$ connecting $z_{\mathfrak a}$ to $w_{\mathfrak a}$ and containing no arrow of type $\mathfrak a$. 
Then $\ga_{\mathfrak a}\mu_{\mathfrak a}$ is a path connecting $z_{\mathfrak a}\in H$ to $v$ and containing no arrow of type $\mathfrak a$. 
As this holds for each $\mathfrak a\in T$, it follows that $v\in P(H)$. As this holds for each $v\in P(P(H))$, we have $P(P(H))\subseteq P(H)$.
\end{proof}

There are other ways of characterizing $P(H)$:

\begin{proposition}\label{PH123} Let $H$ be a non-empty set of vertices, and $v$ a vertex of a $\mathbb{Z}^n$-quiver $Q$. The following statements are equivalent:
  \begin{enumerate}
  \item $v\not\in P(H)$.
  \item There is a vertex $w$ of $Q$ different from $v$ such that for each $z\in H$ there is an admissible path from $z$ to $v$ passing through $w$.
    \item There is a vertex $w$ of $P(H)$ different from $v$ such that for each $z\in H$ there is an admissible path from $z$ to $v$ passing through $w$. (If $P(H)=H$ then $w$ is unique.)
  \item There are a vertex $w$ of $Q$ and a non-empty proper set $I$ of types such that $v\in C_I(w)-\{w\}$ and each admissible path $\gamma$ connecting a vertex of $H$ to $w$ contains no arrow of type $\mathfrak a$ for some $\mathfrak a\not\in I$. \end{enumerate}
  \end{proposition}

  \begin{proof} Assume Statement~(1). By definition, there is an arrow type such that for each $z\in H$ each admissible path $\ga$ connecting $z$ to $v$ contains an arrow of that type. Let $a$ be the arrow of that type arriving at $v$ and $w$ its initial vertex. Given $z$ and $\ga$ as above, let $\mu$ be a path arriving at $w$ such that $a\mu$ has the same type as $\ga$. Since $a \mu$ and $\ga$ arrive at the same vertex, $v$, they leave from the same vertex, $z$. So $a \mu$ is an admissible path from $z$ to $v$ passing through $w$. Statement~(2) holds.

Assume Statement (2). If $w\not\in P(H)$ then Statement (2) applies to $w$ as well, in place of $v$. We may apply it repeatedly. As $H$ is non-empty, and all admissible paths from a given vertex of $H$ to $v$ have the same finite length, the argument cannot be repeated indefinitely. Thus there is a vertex $w\in P(H)$ such that for each $z\in H$ there is an admissible path from $z$ to $v$ passing through $w$. Statement~(3) holds. Furthermore, if $w'$ is another vertex of $P(H)$ with the same property, and $P(H)=H$, then there is an admissible path from $w$ to $v$ passing through $w'$ and one from $w'$ to $v$ passing through $w$. Then an admissible path from $w'$ to $v$ has length at most and at least that of an admissible path from $w$ to $v$, with equality only if $w'=w$, respectively. It follows that $w'=w$. 

Assume Statement~(3). Let $\mu$ be an admissible path connecting $w$ to $v$ and let $I$ be its essential type. Then $I$ is a non-empty proper set of types and $v\in C_I(w)-\{w\}$. Now, by Statement~(3), for each $z\in H$ there are paths $\nu_1$ connecting $z$ to $w$ and $\nu_2$ connecting $w$ to $v$ such that $\nu_2\nu_1$ is admissible. But $\nu_2$ is of the same type as $\mu$. Hence there must be a type $\mathfrak a$ not in $I$ which $\nu_1$ does not contain. This must be true of any admissible path connecting $z$ to $w$. Statement~(4) holds.

Finally, assume Statement~(4). Let $\mu$ be an admissible path connecting $w$ to $v$ and $J$ its essential type. Since $v\in C_I(w)-\{w\}$ we have that $J$ is non-empty and $J\subseteq I$. For each $z\in H$ let $\ga$ be an admissible path from $z$ to $w$. By Statement~(4), there is a type $\mathfrak a$ not in $I$ which $\ga$ does not contain. It follows that $\mu\ga$ is an admissible path connecting $z$ to $v$ and containing all the arrow types in $J$. As this is true for each $z\in H$, we have that $v\not\in P(H)$.
\end{proof}

Let $Q$ be a $\mathbb Z^n$-quiver. Given a vertex $v$ of $Q$ and a collection $I$ of arrow types, all simple paths $\nu$ leaving $v$ and having essential type $I$ end up at the same vertex $w$. We denote $I\cdot v:=w$.

\begin{Def}\label{defpolygon} A finite nonempty collection of pairwise neighboring vertices of a 
$\mathbb Z^n$-quiver is called a \emph{polygon}. If it consists of $m$ vertices, it is called a \emph{$m$-gon}.
\end{Def}

A vertex is a $1$-gon. Neighboring vertices form a $2$-gon, or a \emph{segment}. A 3-gon is also called a \emph{triangle}. 

One way of constructing a $m$-gon for each $m=2,\dots,n+1$ containing a given vertex $v$ of $Q$ is to 
pick a nontrivial ordered partition $T=I_1\cup\dots\cup I_m$ of the set of all arrow types $T$ and put 
$v_1:=v$ and $v_{i+1}=I_i\cdot v_i$ for $i=1,\dots,m-1$. It will follow that $v_1=I_m\cdot v_m$ and that $\{v_1,\dots,v_m\}$ is a $m$-gon. We say that $v_1,\dots,v_m$ \emph{form an oriented $m$-gon}. It is actually the only way:

\begin{proposition}\label{orientedmgon} Let $\Delta$ be a $m$-gon for $m\geq 2$ and $v\in\Delta$. Then there is a unique ordering $v_1,\dots,v_m$ of the vertices of $\Delta$ such that $v_1=v$ and that $v_1,\dots,v_m$ form an oriented $m$-gon. In particular, $m\leq n+1$ and  every $p$-gon is contained in a $(p+1)$-gon for each $p=1,\dots,n$.
\end{proposition}

\begin{proof} We prove the existence of the ordering first. Consider an ordering $v_1,\dots,v_m$ of the vertices of $\Delta$ with $v_1=v$. Since $v_i$ is a neighbor to $v_1$ for each $i>1$, there is a nonempty proper collection of arrow types $J_i$ such that $v_i=J_i\cdot v_1$. Let $I_1:=J_i$ for any $i$ such that $\#J_i=\min\{\# J_j\}$.  Reordering, we may assume $i=2$, whence $v_2=I_1\cdot v_1$.  Suppose by induction that there are nonempty pairwise disjoint proper collections $I_1,\dots,I_p$ of arrow types with $1\leq p < m$ such that $v_{i+1}=I_i\cdot v_i$ for $i=1,\dots,p$ and such that for each $i=1,\dots,p$ and each $j>i$, the collection of arrow types $J$ satisfying $v_j=J\cdot v_i$ satisfies $\# J\geq \# I_i$. 

If $p=m-1$ then $v_1=I_m\cdot v_m$, where $I_m$ is the collection of arrow types not in $I_i$ for any $i$, finishing the proof. If $p<m-1$, since $v_i$ is a neighbor to $v_{p+1}$ for each $i>p+1$, there is a nonempty proper collection of arrow types $L_i$ such that $v_i=L_i\cdot v_{p+1}$. Let $I_{p+1}:=L_i$ for any $i$ such that $\#L_i=\min\{\# L_j\}$.  Reordering, we may assume $i=p+2$, whence $v_{p+2}=I_{p+1}\cdot v_{p+1}$. To finish the proof we need only show now that $I_{p+1}$ is disjoint from $I_i$ for all $i\leq p$.

Now, for each $i=1,\dots,p$, since 
$v_{i+1}$ and $v_{p+2}$ are neighbors, there is an admissible path from $v_i$ to $v_{p+2}$ through $v_{i+1}$ or an admissible path from $v_i$ to $v_{i+1}$ through $v_{p+2}$, by Lemma~\ref{adjadm}. The latter cannot occur, as it would otherwise imply that $v_{p+2}$ is ``closer" to $v_i$ than $v_{i+1}$. So there is an admissible path from $v_i$ to $v_{p+2}$ through $v_{i+1}$ for each $i=1,\dots,p$. For $i=p$ this implies that $I_{p+1}$ is disjoint from $I_p$. Suppose by descending induction on 
$m$ with $1<m\leq p$ that $I_{p+1}$ is disjoint from $I_i$ for  all $i=m,\dots,p$. It follows that $v_{p+2}=(I_m\cup\cdots\cup I_{p+1})\cdot v_m$. Since there is an admissible path from $v_{m-1}$ to $v_{p+2}$ through $v_m$ we have that $I_{m-1}\cap(I_m\cup\cdots\cup I_{p+1})=\emptyset$, and in particular that $I_{p+1}$ is disjoint from $I_{m-1}$.

As for uniqueness, it is enough to observe that, if $v_1,\dots,v_m$ is an ordering of the vertices of $\Delta$ as stated, then $v_{i+1}$ is the closest vertex of $\Delta-\{v_i\}$ to $v_i$ for each $i=1,\dots,m-1$.
\end{proof}

Given vertices $v_1,\dots,v_m$ forming an oriented polygon, and a permutation $\sigma$ of $\{1,\dots,m\}$, it follows from the definition that $v_{\sigma(1)},\dots,v_{\sigma(m)}$ form an oriented polygon if and only if $\sigma$ is a power of $(1 2 \cdots m)$. A subsequence of $v_{\sigma(1)},\dots,v_{\sigma(m)}$ for some power $\sigma$ of $(1 2 \cdots m)$ is said to \emph{form an oriented polygon compatible with that formed by $v_1,\dots,v_m$}.

\begin{proposition}\label{hullpolygon} A polygon is equal to its hull.
\end{proposition}

\begin{proof} Let $v_1,\dots,v_m$ be vertices forming an oriented polygon and $v$ a distinct vertex. As usual, put $v_{m+1}:=v_1$. Let $i$ be such that the length of an admissible path $\gamma$ connecting $v_i$ to $v$ is minimum. Let $I$ be the collection of arrow types 
such that $v_{i+1}=I\cdot v_i$. By the minimality of $i$, there is $\mathfrak a\in I$ not contained in the essential type of $\gamma$. It follows that there is an admissible path from $v_j$ to $v$ passing through $v_i$ for every $j$. Then Proposition~\ref{PH123} yields that $v$ is not in the hull of the polygon.
\end{proof}

\section{Properties of linked nets}\label{generalnets}

%A simple basis of a linked net of vector spaces is just a collection of vectors distributed in the vector spaces of the linked net such that these vectors, when put together, can generate any vector spaces of the linked net. Hence, the existence of a simple basis for a linked net $\g$ gives a decomposition of $\g$ as a direct sum of subnets of pure dimension one. Which simplifying not only the study of $\g$ but also the study of its \emph{linked projective space}.

\begin{Lem}\label{nonzeroadm} Let $\g$ be a weakly linked net over a $\mathbb{Z}^n$-quiver $Q$. If $\gamma$ is a non-admissible path in $Q$ then $\vp_\ga =0$. In particular, 
$\vp^v_u\vp^u_v=0$ for each two distinct vertices $u$ and $v$ of $Q$. 
\end{Lem}

\begin{proof} Let $\ga$ be a non-admissible path. Then there is an intermediate path $\mu$ to $\ga$ which is not admissible and has minimal length for this property. We need only prove $\vp_\mu=0$. 

Now, $\mu=\epsilon a$ for a certain arrow $a$. By the minimality of $\mu$, none of the arrows in $\epsilon$ has the same type as $a$, whence $\epsilon$ is admissible. Furthermore, since $\mu$ is not admissible, there is an admissible path $\rho$ connecting the same two vertices as $\epsilon$ such that $\rho=\rho_2\rho_1$ with $\rho_1 a$ a minimal circuit. Then $\vp_{\epsilon}$ is a scalar multiple of $\vp_{\rho}$, and thus $\vp_\mu$ is a scalar multiple of $\vp_\rho\vp_a$. But $\vp_\rho\vp_a=\vp_{\rho_2}\vp_{\rho_1 a}$, which is zero because $\rho_1 a$ is a minimal circuit. Then $\vp_\mu=0$.
\end{proof}

\begin{Lem}\label{vpuvw} Let $\g$ be a weakly linked net over a $\mathbb{Z}^n$-quiver $Q$. For each three vertices $u,v,w$ of $Q$, if $\vp^v_w\vp^u_v\neq 0$ then there is an admissible path from $u$ to $w$ thorugh $v$; in particular, $\vp^v_w\vp^u_v=\vp^u_w$.
\end{Lem}

\begin{proof} Let $\gamma$ (resp.~$\mu$) be an admissible path connecting $u$ to $v$ (resp.~$v$ to $w$). Then $\vp^v_w\vp^u_v$ is the set of nonzero scalar multiples of $\vp_{\mu\gamma}$. If $\vp^v_w\vp^u_v\neq 0$ then $\mu\gamma$ is admissible by Lemma~\ref{nonzeroadm}.
\end{proof}

Let $\g$ be a weakly linked net of vector spaces over a $\mathbb Z^n$-quiver $Q$. If $\g$ is $1$-generated by a finite set $H$, we may remove vertices from $H$, if necessary, to assume that $H$ is minimal. Then $H$ is unique by the next proposition.

\begin{proposition}\label{lm: generating_are_in_support}
  Let $H$ and $H'$ be sets of vertices of $Q$ and $\g$ be a weakly linked net over a $\mathbb{Z}^n$-quiver $Q$ that is $1$-generated by $H$ and $H'$. If $H'$ is minimal then $H'\subseteq H$ and 
  $\vp^v_w$ is not an epimorphism for each two distinct $v,w\in H'$.
\end{proposition}

\begin{proof} Assume $H'$ is minimal. Suppose by contradiction that there exist distinct $v,w\in H'$ such that $\vp^v_w$ is an epimorphism. We claim that $\g$ is $1$-generated by $H'-\{w\}$, a contradiction. Indeed, let $u$ be a vertex of $Q$. By hypothesis, there exists $z\in H'$ such that $\vp^z_u$ is an epimorphism. If $z \neq w$, we are done. If $z = w$ then the composition
  $\vp^w_u\vp^v_w$ is an epimorphism too. It follows from Lemma~\ref{vpuvw} that 
  $\varphi^v_u=\varphi^w_u\varphi^v_w$, whence $\varphi^v_u$ is an epimorphism as well. We may thus remove $w$ from $H'$, contradicting its minimality.

  Now, suppose by contradiction 
  that $v\notin H$ 
  for some $v\in H'$. Since $\g$ is $1$-generated by $H$, there exists $w\in H$ such that $\vp^w_v$ is an epimorphism. It follows from Lemma~\ref{nonzeroadm} that 
  $\vp^v_w$ is zero. Since
  $\g$ is $1$-generated by $H'$ as well, there is $u\in H'$ such that
  $\vp^u_w$ is an epimorphism. Thus $\vp^w_v\vp^u_w$ is an epimorphism, and hence so is $\vp^u_v$ by Lemma~\ref{vpuvw}.  
  But this contradicts what we have just proved above, as $u,v\in H'$. 
\end{proof}

\begin{proposition}\label{simpfin} Let $\g$ be a weakly linked net over a $\mathbb Z^n$-quiver $Q$ and $H$ a set of vertices of $Q$. If 
$\g$ is generated by $H$ then $\g$ is $1$-generated by $P(H)$. In addition, $\g$ is finitely generated if and only if it is $1$-generated by a finite set. 
\end{proposition}

\begin{proof} %If $\mathfrak g$ has support in a finite set of vertices $S$ then $S$ generates $\mathfrak g$: For each vertex $v$ of the quiver,there is $w\in S$ such that $\vp^w_v$ is an isomorphism, and thus $V_v=\text{Im}(\vp^w_v)$. 
Because of Proposition~\ref{HPH}, we need only prove the first statement. Suppose then that $\g$ is generated by $H$. Let
  $v\not\in P(H)$. Then, by Proposition~\ref{PH123}, there is a vertex $w\in P(H)$ such that there is an admissible path from each
  $z\in H$ to $v$ passing through $w$. Thus $\vp_{v}^{z}=\vp^w_v\vp_{w}^{z}$, and in particular 
  $\text{Im}(\vp_v^z)\subseteq\text{Im}(\vp^w_v)$ for each $z\in H$. Since $\g$ is generated by $H$, it follows that $\varphi^w_v$ is an epimorphism.
\end{proof}

\begin{proposition}\label{noseq} Let $\g$ be a weakly linked net over a $\mathbb Z^n$-quiver $Q$. If $\g$ is finitely generated then $\g$ is locally finite.
\end{proposition}

\begin{proof} If $\g$ is finitely generated then it is $1$-generated by a finite set $H$ by Proposition~\ref{simpfin}. Let $v$ be a vertex of $Q$ and $\ell$ be the maximum of the lengths of the admissible paths connecting a vertex of $H$ to $v$. Let $\mu$ be a path arriving at $v$ of length greater than $\ell$. Then $\vp_\mu=0$. Indeed, there is a vertex $u\in H$ such that $\vp^u_w$ is an epimorphism, where $w$ is the inital vertex of $\mu$. Since 
$\mu$ has length greater than $\ell$, there is no admissible path from $u$ to $v$ through $w$, and hence $[\vp_\mu]\vp^u_w=0$ by Lemma~\ref{vpuvw}. Since 
$\vp^u_w$ is an epimorphism, it follows that $\vp_\mu=0$.
%Let $v$ be a vertex of $Q$ and $\mu_1,\mu_2,\dots$ a sequence of paths arriving at $v$ with increasing lengths. If $\mu_i$ is not admissible then $\vp_{\mu_i}=0$ by Lemma~\ref{nonzeroadm}. Assume all $\mu_i$ are admissible. For each $i$ let $v_i$ be the initial vertex of $\mu_i$. If $\g$ is finitely generated then it is $1$-generated by a finite set by Proposition~\ref{simpfin}, and hence there are a vertex $u$ of $Q$ and an infinite subsequence $v_{i_\ell}$ such that $\vp^u_{v_{i_\ell}}$ is an epimorphism for each $\ell$. We may assume the subsequence is the full sequence. Were $\vp_{\mu_i}$ nonzero for every $i$, then $\vp^{v_i}_{v}\vp^u_{v_i}$ would be nonzero and hence there would exist an admissible path from $u$ to $v$ through $v_i$ for every $i$. This is not possible as all admissible paths from $u$ to $v$ have the same length and the length of an admissible path from $v_i$ to $v$ grows to infinity as $i$ grows.
\end{proof}

\begin{Lem}%[{Santos, \cite[lemma 3.2]{Renan2021}}]
\label{lem_renan} 
Let $\g$ be a weakly linked net over a $\Z^n$-quiver $Q$. Then the following two statements are equivalent:
\begin{enumerate}
    \item $\g$ is a linked net.
    \item $\text{Ker}(\varphi_{\nu_1})\cap\text{Ker}(\varphi_{\nu_2})=0$ for each two simple admissible paths $\nu_1$ and $\nu_2$ in $Q$ with no arrow type in common leaving the same vertex.
    \end{enumerate}
They imply: 
    \begin{enumerate}
    \item[(3)] For each vertex $v$ and arrows $a$ and $b$ of $Q$ of the same type leaving $v$ and arriving at $v$, respectively, 
    $$
    \text{Ker}(\varphi_a)\cap\text{Im}(\varphi_b)=0,
    $$
    \end{enumerate}
which in turn implies:
    \begin{enumerate}
    \item[(4)] For each two admissible paths $\ga_1$ and $\ga_2$ in $Q$ of the same essential type leaving the same vertex,
 \begin{equation}\label{KK}
 \text{Ker}(\vp_{\ga_1})=\text{Ker}(\vp_{\ga_2}).
 \end{equation}
\end{enumerate}
Furthermore, if $\g$ is exact then all four statements are equivalent.
\end{Lem}

\begin{proof} Statement~1 clearly implies Statement~2, which clearly implies Statement~3, as $\varphi_\nu\varphi_b=0$ for each path $\nu$ such that $\nu b$ is a minimal circuit. Also, Statements~2~and~4 clearly imply Statement 1. 
We will now show Statement~3 implies Statement~4.

Assume Statement 3. We make a claim. Let $\mu$ be a nontrivial admissible path. Let $w$ be its final vertex and $a$ an arrow leaving $w$ of the same type of an arrow in $\mu$. We claim that
$$
\text{Ker}(\vp_a) \cap \Ima(\vp_\mu) =0.
$$
Indeed, let $b$ the arrow of the same type as $a$ arriving at $w$. Let $v$ be its initial vertex. Since there is an arrow of the same type as $b$ in $\mu$, there is an admissible path $\beta$ 
such that $b\beta$ is admissible and connects the same two vertices as $\mu$. But then $\text{Im}(\varphi_\mu)=
\text{Im}(\varphi_b\varphi_\beta)\subseteq\text{Im}(\vp_b)$. We may thus apply Statement~3 to prove the claim.

We prove Statement~4. Let $\ga_1$ and $\ga_2$ be admissible paths in $Q$ of the same essential type leaving a vertex $v$. We may assume $\gamma_2$ is simple. We proceed by induction on the length of $\gamma_1$. If $\gamma_1$ has length 0 or 1, Statement 4 holds trivially as then $\ga_1=\ga_2$. It holds as well if $\gamma_1$ is simple, as then $\vp_{\ga_1}$ is a nonzero scalar multiple of $\vp_{\ga_2}$. So we may assume there is $\mathfrak a\in T$ such that $\mathfrak{t}_{\ga_1}(\mathfrak a)\geq 2$. Now, there is an admissible path $\beta$ leaving $v$ whose last arrow $b$ is of type $\mathfrak a$ such that $\mathfrak{t}_{\ga_1}(\mathfrak b)=\mathfrak{t}_{\beta}(\mathfrak b)$ for each $\mathfrak b\neq \mathfrak a$ and 
$\mathfrak{t}_{\ga_1}(\mathfrak a)=\mathfrak{t}_{\beta}(\mathfrak a)+1$. In particular, $\beta$ has the same essential type as $\ga_1$ and thus as $\ga_2$. By induction, $\text{Ker}(\vp_\beta)=\text{Ker}(\vp_{\ga_2})$. 
Let now $a$ be an arrow of type $\mathfrak a$ leaving the end vertex of $b$. Then $\varphi_{\ga_1}$ is a nonzero scalar multiple of $\vp_a\vp_\beta$, whence $\text{Ker}(\vp_{\ga_1})=\text{Ker}(\vp_a\vp_\beta)$. By Statement~3, $\text{Ker}(\vp_a\vp_\beta)=\text{Ker}(\vp_\beta)$. Statement 4 is proved.

Assume now that $\g$ is exact and Statement~4 holds. Let $\beta_1$ and $\beta_2$ be two nontrivial admissible paths in $Q$ leaving the same vertex $v$ of disjoint essential types. Let $\nu_2$ be a simple path of the same essential type as $\beta_2$ and $\mu_2$ a reverse path to $\nu_2$. As $\g$ is exact and Statement~4 holds, $\text{Ker}(\vp_{\beta_2})=\text{Ker}(\vp_{\nu_2})=\text{Im}(\vp_{\mu_2})$. Thus 
$$
\text{Ker}(\vp_{\beta_1})\cap\text{Ker}(\vp_{\beta_2})=\vp_{\mu_2}\big(
\text{Ker}(\vp_{\beta_1}\vp_{\mu_2})\big).
$$
Now, the essential type of $\beta_1$ is contained in the essential type of $\mu_2$, whence 
$\beta_1\mu_2$ is admissible of the same essential type as $\mu_2$. Thus, Statement~4 implies that
$$
\text{Ker}(\vp_{\beta_1}\vp_{\mu_2})=
\text{Ker}(\vp_{\beta_1\mu_2})=\text{Ker}(\vp_{\mu_2}),
$$
and hence $\text{Ker}(\vp_{\beta_1})\cap\text{Ker}(\vp_{\beta_2})=0$. Statement~1 is proved.
\end{proof}

\begin{Def} A subrepresentation of a (weakly) linked net is called a (weakly) linked \emph{subnet}.
\end{Def}

\begin{proposition}\label{subnet} Let $\g$ and $\g'$ be weakly linked net over a $\mathbb Z^n$-quiver and $\hh$ a subnet of $\g$. Then the following statements hold:
\begin{enumerate}
%    \item $\g^*$ is a weakly linked net.
%    \item If $\g$ is special (resp.~general, resp.~exact), so is $\g^*$.
    \item $\g$ and $\g'$ are exact if and only if $\g\oplus\g'$ is.  
    \item $\g$ and $\g'$ are linked nets if and only if $\g\oplus\g'$ is.  
    \item $\hh$ and $\g/\hh$ are weakly linked nets.
    \item If $\g$ is a linked net so is $\hh$.
    \item If $\g$ is locally finite so are $\hh$ and $\g/\hh$.
%    \item If $\g$ is special (resp.~general) so are $\hh$ and $\g/\hh$.
    \item If $\hh$ and $\g/\hh$ are generated by a collection of vertices $H$, so is $\g$.
    \item If $\g$ is $m$-generated by a collection of vertices $H$, so is $\g/\hh$; and so is $\hh$ if $m=1$ and $\g/\hh$ is pure.
    \item If $\g$ and $\hh$ are exact, so is $\g/\hh$.
\end{enumerate}
\end{proposition}

\begin{proof} The first five statements are clearly true. The sixth follows directly from the Snake Lemma. The first statement in (7) is clear, whereas the second follows directly from the Snake Lemma as well. 

Finally, the last statement follows from the Snake Lemma. Indeed, 
assume that $\g$ and $\hh$ are exact. 
Let $v_1$ and $v_2$ be neighboring vertices of $Q$. Let $\mu$ be a simple path connecting $v_1$ to $v_2$ and $\nu$ a reverse path. We need prove that 
$\text{Ker}(\vp^{\mathfrak V/\mathfrak W}_\nu)=\text{Im}(\vp^{\mathfrak V/\mathfrak W}_\mu)$. 
By the Snake Lemma, we have the following diagram of vertical exact sequences:
$$
\begin{CD}
@. @. @. \text{Ker}(\vp^{\mathfrak V}_\mu) @.\\
@. @. @. @VVV @.\\
@. 0 @. 0 @. \text{Ker}(\vp^{\mathfrak V/\mathfrak W}_\mu) @.\\
@. @VVV @VVV @VVV @.\\
0 @>>> \text{Ker}(\vp^{\mathfrak W}_\nu) @>>> V^{\mathfrak W}_{v_2} @>>> 
\text{Coker}(\vp^{\mathfrak W}_\mu) @>>> 0\\
@. @VfVV @VVV @Vg_\mu VV @.\\
0 @>>> \text{Ker}(\vp^{\mathfrak V}_\nu) @>>> V^{\mathfrak V}_{v_2} @>>> 
\text{Coker}(\vp^{\mathfrak V}_\mu) @>>> 0\\
@. @VhVV @VVV @VVV @.\\
0 @>>> \text{Ker}(\vp^{\mathfrak V/\mathfrak W}_\nu) @>>> V^{\mathfrak V/\mathfrak W}_{v_2} @>>> 
\text{Coker}(\vp^{\mathfrak V/\mathfrak W}_\mu) @>>> 0\\
@. @VVV @VVV @VVV @.\\
@. \text{Coker}(\vp^{\mathfrak W}_\nu) @. 0 @. 0 @.\\
@. @Vg_\nu VV @. @. @.\\
@. \text{Coker}(\vp^{\mathfrak V}_\nu) @. @. @.
\end{CD}
$$
Since $\g$ and $\hh$ are exact, the topmost two horizontal short sequences are exact. The bottom one is exact but, arguably, at the middle. Then the diagram is a diagram of exact sequences but, arguably, at $V^{\mathfrak V/\mathfrak W}_{v_2}$. We need to prove the diagram is exact there as well. Again by the Snake Lemma, we have the following long exact sequence:
$$
0\to\text{Ker}(g_\mu) \to \text{Coker}(f) \to V^{\mathfrak V/\mathfrak W}_{v_2} \to 
\text{Coker}(\vp^{\mathfrak V/\mathfrak W}_\mu) \to 0.
$$
But what we know so far of the exactness of the diagram yields that the middle map above is a monomorphism, and thus $\text{Ker}(g_\mu)=0$. Exchanging $\mu$ and $\nu$, by analogy, $\text{Ker}(g_\nu)=0$, and hence $h$ is an epimorphism. Thus 
$\text{Coker}(f)=\text{Ker}(\vp^{\mathfrak V/\mathfrak W}_\nu)$. But then the above long exact sequence is essentially the bottom horizontal sequence of the diagram, which is thus exact.
\end{proof}

\begin{Def} A subrepresentation of a (weakly) linked net is called a (weakly) linked \emph{subnet}.
\end{Def}

The quotient of a pure linked net of vector spaces by a pure subnet need not be a linked net; see Proposition~\ref{subrep_quorep} below though.

\section{Binary linked nets}\label{sec_SB}

%\begin{proposition}\label{subgen} Let $\g$ be a weakly linked net of vector spaces over a $\mathbb{Z}^n$-quiver $Q$. Given vertices $w_1,\dots,w_m$ of $Q$ and sections $s_i\in V_{w_i}$ for $i=1,\dots,m$, for each vertex $v$ of $Q$ let $W_v\subseteq V_v$ be the subspace generated by $\vp_{v}^{w_1}[s_1]\cup\cdots\cup\vp_v^{w_m}[s_m]$. Then the $W_v$ form a subnet $\hh$ of $\g$.
%\end{proposition}

%\begin{proof} It is enough to observe that $\vp^v_u\vp_{v}^{w_j}[s_j]$ is either zero or equal to $\vp_{u}^{w_j}[s_j]$ for each two vertices $u$ and $v$ of $Q$ and $j=1,\dots,m$.
%\end{proof}

\begin{proposition}\label{subgen} Let $\g$ be a weakly linked net over a $\mathbb{Z}^n$-quiver $Q$. Given vertices $v_1,\dots,v_m$ of $Q$ and subojects  $W_i\subseteq V_{v_i}$ for $i=1,\dots,m$, for each vertex $v$ of $Q$ let $W_v:=\sum\vp_v^{w_i}(W_i)\subseteq V_v$. Then the $W_v$ form a subnet $\hh$ of $\g$.
\end{proposition}

\begin{proof} It is enough to observe that $\vp^v_u\vp_{v}^{w_j}$ is either zero or 
equal to $\vp_{u}^{w_j}$ for each two vertices $u$ and $v$ of $Q$ and $j=1,\dots,m$.
\end{proof}

%\begin{Def}\label{def_SB} Given the setup of Proposition~\ref{subgen}, we say that the collection of vertices $v_i$ of $Q$ and subojects $W_i\subseteq V_{v_i}$ \emph{generate} $\hh$. If $m=1$ we say that $\hh$ is \emph{unigenerated}, and if in addition the epimorphism $W_{w_1}\to W_v$ is an isomorphism for each vertex $v$ of $Q$ we say that $\hh$ is \emph{faithfully unigenerated}. 
%We will also say that $\{v_1,\dots,v_m\}$ \emph{generate} $\mathfrak g$ and that $\mathfrak g$ is finitely generated. 
%\end{Def}

%\begin{Def}\label{def_SB} Given the setup of Proposition~\ref{subgen}, we say that the collection of vertices $w_1,\dots,w_m$ of $Q$ and vectors $s_i\in V_{w_i}$ for $i=1,\dots,m$ \emph{generate} $\hh$. 
%We will also say that $\{v_1,\dots,v_m\}$ \emph{generate} $\mathfrak g$ and that $\mathfrak g$ is finitely generated. If in addition $\hh$ is pure of dimension $m$ then we say that they form a \emph{simple basis} for $\hh$, and that $\hh$ admits a simple basis.
%\end{Def}

%Of course, if $\g$ admits a (finite) set of generators then $\g$ has finite dimension, and if $\g$ admits a simple basis then $\g$ has finite pure dimension.

\begin{Def}\label{def_binary} A weakly linked net $\g$ is \emph{binary} if each of its associated maps is either zero or a monomorphism. It is called \emph{simple} if each of its associated objects is simple.
\end{Def}

A simple weakly linked net is clearly pure and binary. A weakly linked net of vector spaces is simple if and only if all the associated spaces have dimension 1.

\begin{proposition}\label{binexact} A weakly linked net $\g$ is exact and binary if and only if for each two neighboring vertices $u$ and $v$, either $\vp^u_v$ or $\vp^v_u$ is an isomorphism. 
\end{proposition}

\begin{proof} Let $u$ and $v$ be neighboring vertices. Since $\vp^v_u\vp^u_v=0$, if $\vp^u_v$ is an isomorphism, then $\vp^v_u=0$ and hence $\text{Im}(\vp^u_v)=V_v=\text{Ker}(\vp^v_u)$. On the other hand, if $\vp^v_u$ is an isomorphism then $\vp^u_v=0$, and hence $\text{Im}(\vp^u_v)=0=\text{Ker}(\vp^v_u)$. Either way, $\g$ is binary and exact. 

Conversely, suppose $\g$ is binary and exact. Then $\text{Im}(\vp^u_v)=\text{Ker}(\vp^v_u)$ and 
$\text{Im}(\vp^v_u)=\text{Ker}(\vp^u_v)$. It follows that $\vp^v_u$ is an isomorphism if $\vp^u_v=0$. On the other hand, if $\vp^u_v\neq 0$ then $\vp^u_v$ is a monomorphism because $\g$ is binary. Since $\vp^u_v\vp^v_u=0$, we have $\vp^v_u=0$. Mutatis mutandi, $\vp^u_v$ is an isomorphism.
\end{proof}

%This Lemma implies that if $\vp^{v_1}_{v_2}$ is an isomorphism, then $\vp^{v_2}_{v_3}$ is also an isomorphism, and if $\vp^{v_2}_{v_3}$ is the null map then $\vp^{v_1}_{v_2}$ is null as well.

\begin{Def} Let $\g$ be a weakly linked net over a $\mathbb{Z}^n$-quiver $Q$, and $v$ a vertex of $Q$. We say $v$ is \emph{primitive} for $\g$ if 
$$
V_v\neq\sum_a\text{\rm Im}(\varphi_a),
$$
where $a$ runs through all arrows arriving at $v$.
\end{Def}

\begin{Lem}\label{subrep_prim} Let $\g$ be a locally finite weakly linked net over a $\mathbb Z^n$-quiver $Q$ and $v$ a vertex of $Q$. If $V_v\neq 0$ then there are primitive vertices $v_1,\dots,v_m$ of $Q$ for $\g$ such that $V_v=\sum\text{Im}(\vp^{v_i}_v)$.
\end{Lem} 

\begin{proof} Suppose by contradiction that the statement does not hold. Then $v$ is not primitive for $\g$, that is, 
\begin{equation}\label{vvi}
V_v=\sum_i\text{\rm Im}(\varphi_{\mu_i})
\end{equation}
for certain admissible paths $\mu_i$ of length 1 arriving at $v$ and satisfying $\varphi_{\mu_i}\neq 0$.

We claim we can find a sequence of expressions for $V_v$ as in \eqref{vvi}, for certain finitely many admissible paths $\mu_i$ arriving at $v$ with $\vp_{\mu_i}\neq 0$, such that the minimum of the lengths of the $\mu_i$ leaving from nonprimitive vertices for $\g$ goes to infinity, contradicting the hypothesis that $\g$ is locally finite. 

%from the claim: Since $H$ is finite, the lengths of admissible paths connecting vertices of $H$ are bounded. 

To prove the claim, consider an expression as in \eqref{vvi} for finitely many admissible paths $\mu_i$ arriving at $v$ with $\varphi_{\mu_i}\neq 0$ for each $i$. Let $w_i$ be the initial vertex of $\mu_i$ for each $i$. For each $i$ such that $w_i$ is nonprimitive for $\g$, that is, $V_{w_i}=\sum\text{Im}(\varphi_a)$, where the $a$ are the arrows arriving at $w_i$, we have
$$
\text{\rm Im}(\varphi_{\mu_i})=\sum_a\text{\rm Im}(\varphi_{\mu_ia}).
$$
%For each $i$, there are $z_i\in H$, an admissible path $\gamma_i$ connecting $z_i$ to $w_i$ and $t_i\in V_{z_i}$ such that $s_i=\varphi_{\gamma_i}(t_i)$. Since
%$$
%\varphi_{\mu_i}\varphi_{\gamma_i}(t_i)=\varphi_{\mu_i}(s_i)\neq 0,
%$$
%it follows from Lemma~\ref{nonzeroadm} that the concatenation of $\gamma_i$ with $\mu_i$ is admissible. Up to replacing $\mu_i$ with this concatenation, which has at least the same length as $\mu_i$, and $s_i$ with $t_i$ for each $i$, we may assume $w_i\in H$. 
Removing from the above sum those $a$ for which 
$\varphi_{\mu_ia}=0$, in particular, those $a$ such that $\mu_ia$ is not admissible, and substituting in \eqref{vvi} for each $i$, we obtain an expression as in \eqref{vvi}. But now the minimum of the lengths of the $\mu_i$ leaving from nonprimitive vertices for $\g$ has increased.
%If one of them does not start in $H$, we replace it by one that does, as explained above. The claim is proved.
\end{proof}

\begin{proposition} \label{simpleimpliesexact} A weakly linked net is faithfully generated by one vertex if and only if it is a locally finite binary exact linked net.
%or its maps associated to minimal circuits are nonzero multiple of the identity map. 
\end{proposition} 

\begin{proof} 
%If the maps associated to minimal circuits are nonzero there is nothing to prove.
%
%If there is a minimal circuit $\gamma$ with $\vp_\gamma=0$, the $\vp_\mu=0$ for all circuit $\mu$ by Lemma~\ref{lemma_zero_circuits}. In this case
Let $\g$ be a weakly linked net over a $\mathbb{Z}^n$-quiver $Q$. Assume $\g$ is faithfully generated by a vertex $w$ of $Q$. Clearly, $\g$ is finitely generated, whence locally finite by Proposition~\ref{noseq}.

For each two vertices $u$ and $v$, we have $\varphi^u_v \vp_{u}^{w} = \vp_{v}^{w}$ if there is an admissible path from $w$ to $v$ passing through $u$. In this case, since $\vp_{u}^w$ and $\vp_v^w$ are isomorphisms, so is $\vp^u_v$. By Lemma~\ref{adjadm}, 
if $u$ and $v$ are neighbors, if there is no admissible path from $w$ to $v$ through $u$, there is one from $w$ to $u$ through $v$. Thus either $\vp^u_v$ or $\vp^v_u$ is an isomorphism. It follows from Proposition~\ref{binexact} that $\g$ is binary and exact.

Finally, given two admissible paths $\gamma_1$ and $\gamma_2$ leaving the same vertex $u$, the kernels of 
$\varphi_{\gamma_1}$ and $\varphi_{\gamma_2}$ intersect nontrivially if and only if there is an admissible path $\nu$ 
connecting $w$ to $u$ such that  
$\gamma_1\nu$ and $\gamma_2\nu$ are not admissible. But since $\nu$ is admissible, there is an arrow type $\mathfrak a$ such that $\mathfrak{t}_{\nu}(\mathfrak a)=0$, and hence $\mathfrak{t}_{\gamma_1}(\mathfrak a)>0$ and $\mathfrak{t}_{\gamma_2}(\mathfrak a)>0$. The third property of a linked net is verified. 

Conversely, let $\g$ be a locally finite binary exact linked net. If $\g$ is trivial then $\g$ is clearly generated by any vertex. Assume $\g$ is nontrivial. Then $\g$ admits a primitive vertex $v$ by Lemma~\ref{subrep_prim}. Then, for each arrow $a$ arriving at $v$, the morphism $\vp_a$ is not an epimorphism. 
Since $\g$ is binary and exact, it follows from Proposition~\ref{binexact} that $\vp_{\mu_a}$ is an isomorphism, where $\mu_a$ is a reverse path to $a$. Since $\g$ is a linked net, it follows from Lemma~\ref{lem_renan} that $\vp^v_u$ is a monomorphism for each vertex $u$ of $Q$. 

It remains to prove that $\vp^v_u$ is an epimorphism for each vertex $u$ of $Q$. Suppose not. Let $u$ be a vertex of $Q$ for which $\vp^v_u$ is not an epimorphism. Then there is a primitive vertex $w$ for $\g$ distinct from $v$ by Lemma~\ref{subrep_prim}. Mutatis mutandi, 
$\vp^w_v$ is a monomorphism. Since $\vp^v_w$ is a monomorphism as well, it follows that $V_w=0$. But then $w$ is not primitive, a contradiction.
\end{proof}

\begin{Lem}\label{subrep_dim_1} 
Let $\g$ be a linked net over a $\Z^n$-quiver $Q$. Let $v$ be a vertex of $Q$ and $W\subseteq V_v$ a subobject. Let $\hh$ be the subnet generated by $W$. Then $\hh$ is faithfully generated by $v$ only if $W\cap \text{\rm Im}(\varphi_a)=0$ for each arrow $a$ arriving at $v$. The converse holds if $\g$ is exact.
\end{Lem} 

\begin{proof} Of course, $\hh$ is generated by $v$. If $W\cap \text{\rm Im}(\varphi_a)\neq 0$ for an arrow $a$ arriving at $v$ then $\vp_{\mu_a}|_W$ is not a monomorphism for any reverse path $\mu_a$ to $a$, and thus $\hh$ is not faithfully generated by $v$.

On the other hand, assume $\g$ is exact and $W\cap \text{\rm Im}(\varphi_a)=0$ for each arrow $a$ arriving at $v$. We need prove that $\vp^v_w|_W$ is a monomorphism for each vertex $w$ of $Q$. But $\text{Ker}(\vp^v_w)=\text{Ker}(\vp^v_u)$ for a neighbor $u$ of $v$ by Lemma~\ref{lem_renan}, since $\g$ is a linked net. Since $\g$ is exact, $\text{Ker}(\vp^v_u)=\text{Im}(\vp^u_v)$. Now, $\text{Im}(\vp^u_v)\subseteq\text{Im}(\vp_a)$ for an arrow $a$ arriving at $v$. Thus $W\cap\text{Ker}(\vp^v_u)=0$ by hypothesis. 
%
%Furthermore, $\mathfrak h$ is exact. Indeed, if $w_1$ and $w_2$ are adjacent vertices, then there is an admissible path connecting $v$ to $w_1$ through $w_2$ if and only if there is no admissible path connecting $v$ to $w_2$ through $w_1$. Thus, if $\varphi^{w_1}_{w_2}\varphi^{v}_{w_1}(s)=0$ then the latter holds, and the former implies $\varphi^v_{w_1}(s)=\varphi^{w_2}_{w_1}\varphi^{v}_{w_1}(s)$.
\end{proof}

\begin{theorem}\label{simple1} 
A simple locally finite exact linked net is faithfully generated by a vertex.
\end{theorem}

\begin{proof} Let $\g$ be the linked net of the statement. Then $\g$ is nontrivial. By Lemma~\ref{subrep_prim}, there is a primitive vertex $v$ for $\g$. Since $V^\g_v$ is simple, 
$\text{Im}(\vp^\g_a)=0$ for each arrow $a$ arriving at $v$. Applying 
Lemma~\ref{subrep_dim_1} with $W=V^\g_v$ we obtain a subnet $\hh$ of $\g$ with $V^\hh_v=V^\g_v$ which is faithfully generated by $v$. Since $V^\hh_v\neq 0$, also $V^\hh_w\neq 0$ for every vertex $w$ of $Q$. Since $V^\g_w$ is simple, $V^\hh_w=V^\g_w$ for every $w$, that is, $\hh=\g$, and thus 
$\g$ is faithfully generated by $v$.
\end{proof}

%\begin{theorem}\label{simple2} A Noetherian pure exact linked net $\g$ of vector spaces over a $\mathbb Z^n$-quiver admits a simple basis if and only if $\g$ is \emph{completely decomposable}, that is, a direct sum of pure exact linked subnets of dimension 1.
%\end{theorem}

%\begin{proof} The ``only if" statement is clear from the very definition, coupled with Proposition~\ref{simpleimpliesexact} applied to one-dimensional linked nets. The ``if" statement follows from Theorem~\ref{simple1}. 
%\end{proof}

\begin{Cor}\label{simplevec1} A finitely generated exact linked net of vector spaces of dimension~1 is generated by a vertex.
\end{Cor}

\section{The intersection property}
\label{intprop}

Let $Q$ be a $\mathbb Z^n$-quiver. Given a vertex $v$ of $Q$, a collection $I$ of arrow types, and a weakly linked net $\g$ over $Q$, we put $\varphi^v_I:=\varphi^I_{I\cdot v}:=\varphi^v_{I\cdot v}$. If $v$ is clear from the context, we put $\varphi_I:=\varphi^v_I$ and $\varphi^I:=\varphi^I_{I\cdot v}$.

\begin{Lem}\label{ker_ker_ker}
 Let $\g$ be a weakly linked net over a $\Z^n$-quiver. Then $\g$ is a linked net if and only if for each two collections of arrow types $I_1,I_2$ and each vertex $v$ of the quiver,
\begin{equation}\label{I1I2K}
 \text{Ker}(\vp_{I_1}^v)\cap\text{Ker}(\vp_{I_2}^v)=\text{Ker}(\vp_{I_1\cap I_2}^{v}).
\end{equation}
\end{Lem}

\begin{proof} Let $v$ be a vertex of the quiver 
and $I_1,I_2$ collections of arrow types. If $I_i$ is the full set $T$ of arrow types for some $i$, then \eqref{I1I2K} holds clearly. We may now assume $I_1,I_2\subsetneqq T$. Let $w$ be the end vertex of a simple path leaving $v$ with essential type $I_1\cap I_2$. Then $\vp_{I_1}^v=\vp^w_{I_1-I_2}\vp_{I_1\cap I_2}^v$ and 
 $\vp_{I_2}^v=\vp^w_{I_2- I_1}\vp_{I_1\cap I_2}^v$. The inclusion  
 $\text{Ker}(\vp_{I_1}^v)\cap\text{Ker}(\vp_{I_2}^v)\supseteq\text{Ker}(\vp_{I_1\cap I_2}^{v})$ follows.
 
 On the other hand, if $\g$ is a linked net, as $(I_1-I_2) \cap (I_2-I_1)=\emptyset$, we have
 $$
 \text{Ker}(\vp^w_{I_1-I_2}) \cap \text{Ker}(\vp^w_{I_2-I_1})=0.
 $$
 Thus, $\text{Ker}(\vp_{I_1}^v)\cap\text{Ker}(\vp_{I_2}^v)\subseteq\text{Ker}(\vp_{I_1\cap I_2}^v)$ as well, proving \eqref{I1I2K}.
 
Conversely, if \eqref{I1I2K} holds for each two subsets $I_1,I_2\subseteq T$, in particular it holds for proper subsets with empty intersection. Then $\g$ is a linked net by Lemma~\ref{lem_renan}.
\end{proof}

\begin{Def}\label{defintprop} 
A weakly linked net over a $\Z^n$-quiver \emph{satisfies the intersection property} at a vertex $v$ if for all collections $I_0,\dots,I_m$ of arrow types the following 
equality holds:
\begin{equation}\label{int}
\bigg(\sum_{\ell=1}^m\text{Ker}(\varphi_{I_\ell}^v)\bigg)\cap\text{Ker}(\varphi_{I_0}^v) = \sum_{\ell=1}^m\text{Ker}(\varphi_{I_\ell\cap I_0}^v).
\end{equation}
\end{Def}

Notice that \eqref{int} holds trivially if $I_j$ is the full set of arrow types for some $j$.

By Lemma~\ref{ker_ker_ker}, a weakly linked net satisfying the intersection property at every vertex is a linked net. In the light of the lemma, another way of viewing the intersection property for linked nets is by saying that the intersection of kernels distributes with respect to the sum. 

In \cite{G_munoz}, Def.~2.6, p.~556, on limit linear series for curves of compact type with three components, G.~Mu\~noz defines, in a different context, a similar property 
for the special case when $m=2$. 

\begin{Lem}\label{intH} Let $\g$ be a pure linked net over a $\mathbb Z^n$-quiver $1$-generated by a collection of vertices $H$. Then $\g$ satisfies the intersection property at every vertex if $\g$ satisfies it at each vertex of $H$.
\end{Lem}
 
\begin{proof} Let $v$ be a vertex of the quiver not in $H$. Then there is $w\in H$ such that $\varphi^w_v$ is an epimorphism. Let $J$ be the essential type of an admissible path connecting $w$ to $v$. Let 
$I_0,\dots,I_m$ be collections of arrow types. For each $\ell=1,\dots,m$, let $u_\ell$ be the end vertex of a simple path leaving $v$ of essential type $I_\ell$. For each 
$\ell=0,\dots,m$, put $J_\ell:=J\cup I_\ell$. If the linked net satisfies the intersection property at $w$, then
\begin{equation}\label{Jint}
\bigg(\sum_{\ell=1}^m\text{Ker}(\varphi_{J_\ell}^w)\bigg)\cap\text{Ker}(\varphi_{J_0}^w) = \sum_{\ell=1}^m\text{Ker}(\varphi_{J_\ell\cap J_0}^w).
\end{equation}
Now, since $\vp^w_v$ is an epimorphism,
$\text{Ker}(\varphi_{I_\ell}^v)=\vp^w_v(\text{Ker}(\varphi_{u_\ell}^v\varphi^w_v))$ for each $\ell$. If $J_\ell$ is a proper collection of arrow types then 
$\varphi_{u_\ell}^v\varphi^w_v=\varphi^w_{u_\ell}$ and hence
\begin{equation}\label{JJJ}
\text{Ker}(\varphi_{u_\ell}^v\varphi^w_v)=\text{Ker}(\varphi_{J_\ell}^w)
\end{equation}
by Lemma~\ref{lem_renan}; otherwise, $\varphi_{u_\ell}^v\varphi^w_v=0$ and $\vp^w_{J_\ell}=0$, and hence \eqref{JJJ} holds as well. It follows that
$$
\text{Ker}(\varphi_{I_\ell}^v)=\vp^w_v(\text{Ker}(\varphi_{J_\ell}^w))\text{ for each }\ell=1,\dots,m. 
$$
Analogously,
$$
\text{Ker}(\varphi_{I_\ell\cap I_0}^v)=\vp^w_v(\text{Ker}(\varphi_{J_\ell\cap J_0}^w))\text{ for each }\ell=1,\dots,m. 
$$
Finally, since $\g$ is pure, $\vp^w_v$ is an isomorphism, and hence \eqref{int} follows from \eqref{Jint}.
\end{proof}

Given a weakly linked net $\g$ over a $\mathbb Z^n$-quiver $Q$ and a vertex $v$ of $Q$, we call a subobject $W\subseteq V_v$ \emph{primitive} for $\g$ if 
$$
W\cap\sum_{\mathfrak a}\text{Im}\vp^{\mathfrak a}_v=0,
$$
where the sum runs over all arrow types.

\begin{proposition}\label{thm_quo_inter_pro} Let $\g$ be an exact weakly linked net over a $\Z^n$-quiver $Q$ satisfying the intersection property at every vertex. Let $\hh$ be an exact subnet of $\g$ such that for each vertex $v$ of $Q$ there is a vertex $u$ such that $V^\hh_v=\vp^u_v(V^\hh_u)$ and $V^\hh_u$ is primitive for $\g$. Then the quotient $\g/\hh$ is an exact weakly linked net satisfying the intersection property at every vertex. 
\end{proposition}

\begin{proof} That $\g/\hh$ is an exact weakly linked net follows from Proposition~\ref{subnet}. For convenience, we put $W_v:=V^\hh_v$ for each vertex $v$ of $Q$. Let $T$ be the full set of arrow types. Let $v$ be a vertex of $Q$. For each subset $I\subseteq T$, we claim that
\begin{equation}\label{KW}
(\vp^v_{I})^{-1}(W_{I\cdot v})=W_v+\text{Ker}(\vp^v_{I}),
\end{equation}
Indeed, 
$$
(\vp^v_{I})^{-1}(W_{I\cdot v})=(\vp^v_{I})^{-1}\Big(W_{I\cdot v}\cap\text{Ker}(\vp^{I\cdot v}_{T-I})\Big)=(\vp^v_{I})^{-1}\Big(\vp^v_{I}(W_v)\Big)=W_v+\text{Ker}(\vp^v_{I}),
$$
where the middle equality follows from the exactness of $\hh$.

Let $I_0,\dots,I_m$ be proper subsets of $T$. To prove that $\g/\hh$ satisfies the intersection property at $v$ we need to prove that
$$
\bigg(\sum_{\ell=1}^m(\vp^v_{I_\ell})^{-1}(W_{I_\ell\cdot v})\bigg)\cap(\vp^v_{I_0})^{-1}(W_{I_0\cdot v})=\sum_{\ell=1}^m(\vp^v_{I_\ell\cap I_0})^{-1}(W_{(I_\ell\cap I_0)\cdot v}).
$$
By \eqref{KW}, this is equivalent to $(W_v+A)\cap(W_v+B)=W_v+C$, where
$$
A:=\sum_{\ell=1}^m \text{Ker}(\vp^v_{I_\ell}),\quad 
B:=\text{Ker}(\vp^v_{I_0})\quad\text{and}\quad
C:=\sum_{\ell=1}^m\text{Ker}(\vp^v_{I_\ell\cap I_0}).
$$
Since $\g$ is exact, $C=A\cap B$. Hence, we need only show that 
\begin{equation}\label{WABC}
(W_v+A)\cap(W_v+B)=W_v+(A\cap B).
\end{equation}

It is elementary that \eqref{WABC} holds if $W_v\subseteq A$ or $W_v\subseteq B$. 
Suppose $W_v\not\subseteq A$ and $W_v\not\subseteq B$. If we show that 
$W_v\cap(A+B)=0$ then \eqref{WABC} follows. 
Now, $W_v=\vp^u_v(W_u)$ for a vertex $u$ of $Q$ for which $W_u$ is primitive for $\g$. Using that $\g$ is exact at $v$ we get
$$
W_v\cap(A+B)=\vp^u_v\Bigg(
W_u\cap(\vp^u_v)^{-1}\bigg(\sum_{\ell=0}^m\text{Im}(\varphi^{T-I_\ell}_v)\bigg)\Bigg).
$$
Since $\g$ satisfies the intersection property at every vertex, 
$$
W_v\cap(A+B)=\vp^u_v\Bigg(W_u\cap\bigg(\sum_{\ell=0}^m\text{Im}(\varphi^{(T-I_\ell)\cap(T-K)}_u)\bigg)\Bigg)
$$
by Lemma~\ref{SumIma_SumIma} below, where $K$ is the essential type of an admissible path connecting $u$ to $v$. 

Were $(T-I_\ell)\cap(T-K)=\emptyset$ for some $\ell$, then $I_\ell\cup K=T$ and then 
$\vp^v_{I_\ell}(W_v)=\vp^v_{I_\ell}\vp^u_v(W_u)=0$, implying that $W_v\subseteq A$ or $W_v\subseteq B$. As this is not the case, and $W_u$ is primitive for $\g$, we must have 
that $W_v\cap(A+B)=0$.
\end{proof}

%Let $u$ and $v$ be vertices of a $\mathbb Z^n$-quiver of $Q$. Then there are unique vertices $u_0,u_1,\dots,u_p$ such that:
%\begin{enumerate}
%    \item $u_0=u$ and $u_p=v$;
%    \item $u_i=K_i\cdot u_{i-1}$ for a proper collection $K_i$ of arrow types for each $i$.
%    \item $K_1\subseteq %K_2\subseteq\cdots\subseteq K_p$.
%\end{enumerate}
%Indeed, $p:=\max\{t_\mu(\mathfrak a)\}$ for an admissible path $\mu$ connecting $u$ to $v$, and $K_i$ is the essential type of an admissible path connecting $u$ to $u_i$ for $i=p,\cdots,1$.

%\begin{Def} Given two vertices $u$ and $v$ of a $\mathbb Z^n$-quiver of $Q$, we call the sequence $u_0,u_1,\dots,u_p$ above the \emph{lexicographic sequence of intermediate points} from $u$ to $v$.  
%\end{Def}

\begin{Lem}\label{SumIma_SumIma} Let $\g$ be an exact weakly linked net over a $\Z^n$-quiver $Q$. Then $\g$ satisfies the intersection property at every vertex if and only if for each two vertices $u$ and $v$ of $Q$ and all collections $J_1,\dots,J_m$ of arrow types,
\begin{equation}\label{intPequiv}
(\varphi^u_v)^{-1}\Big(\sum_{\ell=1}^m\text{\rm Im}(\varphi^{J_\ell}_v)\Big)=
\sum_{\ell=1}^m\text{\rm Im}(\varphi^{J_\ell\cap(T-K)}_u),
\end{equation}
where $K$ is the essential type of an admissible path connecting $u$ to $v$.
\end{Lem}

\begin{proof} Let $T$ be the collection of all arrow types. Suppose first that $\g$ satisfies the intersection property at every vertex. Let $u$ and $v $ be vertices of $Q$ and $J_1,\dots,J_m\subseteq T$. 
Put $I_\ell:=T-J_\ell$ for $\ell=1,\dots,m$. There are a sequence of vertices $u_0,u_1,\dots,u_p$ and subsets 
$K_1,\dots,K_p\subseteq T$ such that:
\begin{enumerate}
    \item $u_0=u$ and $u_p=v$,
    \item $u_i:=K_i\cdot u_{i-1}$ for $i=1,\dots,p$,
\item $\bigcup K_i=K$.
\end{enumerate}
Then $\vp^u_v=\vp^{u_{p-1}}_{K_p}\cdots\vp^{u_0}_{K_1}$. 

We prove \eqref{intPequiv} by induction on $p$. It is enough to prove it for $p=1$. Then, by the exactness of $\g$,
$$
(\varphi^K_v)^{-1}\Big(\sum_{\ell=1}^m\text{\rm Im}(\varphi^{J_\ell}_v)\Big)=(\varphi^K_v)^{-1}\Bigg(\Big(\sum_{\ell=1}^m\text{Ker}(\varphi^v_{I_\ell})\Big)\cap\text{Ker}(\varphi^v_{T-K})\Bigg),
$$
and thus, by the intersection property at $v$,
$$
(\varphi^K_v)^{-1}\Big(\sum_{\ell=1}^m\text{\rm Im}(\varphi^{J_\ell}_v)\Big)=
(\varphi^K_v)^{-1}\Big(\sum_{\ell=1}^m\text{Ker}(\varphi^v_{I_\ell\cap(T-K)})\Big).
$$
Again, since $\g$ is exact,
$$
(\varphi^K_v)^{-1}\Big(\sum_{\ell=1}^m\text{\rm Im}(\varphi^{J_\ell}_v)\Big)=
(\varphi^K_v)^{-1}\Big(\sum_{\ell=1}^m\text{\rm Im}(\varphi^{J_\ell\cup K}_v)\Big).
$$
But $\varphi^{J_\ell\cup K}_v=\varphi^K_v\varphi^{J_\ell\cap (T-K)}_u$ for each $\ell$. Thus
$$
(\varphi^K_v)^{-1}\Big(\sum_{\ell=1}^m\text{\rm Im}(\varphi^{J_\ell}_v)\Big)=
\text{Ker}(\vp^K_v)+\sum_{\ell=1}^m\text{\rm Im}(\varphi^{J_\ell\cap(T-K)}_u),
$$
whence, again by the exactness of $\g$,
$$
(\varphi^K_v)^{-1}\Big(\sum_{\ell=1}^m\text{\rm Im}(\varphi^{J_\ell}_v)\Big)=
\text{Im}(\vp^{T-K}_u)+\sum_{\ell=1}^m\text{\rm Im}(\varphi^{J_\ell\cap(T-K)}_u).
$$
Since $\varphi^{J_\ell\cap(T-K)}_u$ is a factor of $\vp^{T-K}_u$ for each $\ell$, we obtain \eqref{intPequiv}.

Conversely, let $v$ be a vertex of $Q$ and $I_0,\dots,I_m\subseteq T$. Put $J_\ell:=T-I_\ell$ for $\ell=0,\dots,m$. Put $u:=I_0\cdot v$. Then
\begin{align*}
\bigg(\sum_{\ell=1}^m\text{Ker}(\varphi_{I_\ell}^v)\bigg)\cap\text{Ker}(\varphi_{I_0}^v) & =  \bigg(\sum_{\ell=1}^m\text{Im}(\varphi^{J_\ell}_v)\bigg)\cap\text{Im}(\varphi^{J_0}_v) \\
&=\varphi^{J_0}_v(\varphi^{J_0}_v)^{-1}\bigg(\sum_{\ell=1}^m\text{Im}(\varphi^{J_\ell}_v)\bigg)
= \varphi^{J_0}_v\bigg(\sum_{\ell=1}^m\text{\rm Im}(\varphi^{J_\ell\cap I_0}_u)\bigg) \\
 & =  \sum_{\ell=1}^m\text{\rm Im}(\varphi^{J_\ell\cup J_0}_v)
=\sum_{\ell=1}^m\text{Ker}(\varphi_{I_\ell\cap I_0}^v),
\end{align*}
where exactness of $\g$ is used in the first and last equalities and \eqref{intPequiv} is used in the third equality. The reverse inclusion is trivial. Thus $\g$ satisfies the intersection property at $v$.
\end{proof}

\section[Semisimpleness]{A local criterion for semisimpleness}\label{Classification}

%In this section we prove the main result of the paper, Theorem \ref{cat}, where we characterize with a local condition when a Noetherian exact linked net in a $k$-linear semisimple Abelian category is semisimple.

\begin{proposition}\label{SB_inter}
   A locally finite binary exact linked net satisfies the intersection property at every vertex.
\end{proposition}

\begin{proof} Let $\g$ be the linked net of the statement. By Proposition~\ref{simpleimpliesexact}, there is a vertex $w$ of $Q$ which faithfully generates $\g$. 
Let $v$ be a vertex of $Q$ and $\gamma$ be an admissible path from $w$ to $v$. Let $J$ be the essential type of $\gamma$. Let $T$ be the full collection of arrow types and $I_0,\dots,I_m\subseteq T$. For each subset $I\subseteq T$,
$$
\text{Ker}(\vp_I^{v})=
\begin{cases}
V_v&\text{if }I\cup J=T,\\
0&\text{otherwise.}
\end{cases}
$$
If $I_0\cup J\neq T$ then also $(I_0\cap I_\ell)\cup J\neq T$ for each $\ell$. Thus the intersection property \eqref{int} holds because both sides of the equation are zero. If $I_0\cup J=T$ and $I_\ell\cup J=T$ for some $\ell>0$ then also $(I_\ell\cap I_0)\cup J=T$, and now \eqref{int} holds because both sides of the equation are equal to $V_v$. Finally, if $I_0\cup J=T$ and $I_\ell\cup J\neq T$ for every $\ell>0$, then also $(I_\ell\cap I_0)\cup J\neq T$ for each $\ell>0$, and thus \eqref{int} holds because both sides of the equation are equal to zero.
\end{proof}

\begin{proposition}\label{subrep_quorep}
Let $\g$ be an exact linked net over a $\Z^n$-quiver satisfying the intersection property at every vertex. 
Let $\hh$ be an exact simple subnet of $\g$ such that for each vertex $u$ of $Q$ there is a vertex $v$ such that $V^\hh_u=\vp^v_u(V^\hh_v)$ and $V^\hh_v$ is primitive for $\g$.
%Let $v$ be a vertex of $Q$ and $W\subseteq V_v$ a simple subobject primitive for $\g$. Let $\hh$ be the subnet generated by $W$. 
Then the quotient $\g/\hh$ is an exact linked net.
\end{proposition}

\begin{proof} 
%By Lemma~\ref{subrep_dim_1}, the net $\hh$ is faithfully generated by $v$. By Lemma~\ref{simpleimpliesexact}, the net $\g$ is binary, thus exact by Lemma~\ref{binexact}. It follows now from 
By Proposition~\ref{subnet}, the quotient $\g/\hh$ is exact. We need only prove that it is a linked net. 

For convenience, we put $W_w:=V^\hh_w$ for each vertex $w$ of $Q$. 
%Since $\hh$ is faithfully generated, $W_w=\vp^v_w(W)$ and $W$ is simple, so is $W_w$ for every $w$. 
Let $T$ be the full set of arrow types. Let $u$ be a vertex of $Q$ and $\mathfrak a\in T$. Let $a$ and $b$ be arrows of type $\mathfrak a$ leaving and arriving at $u$, respectively. Let $z$ be the initial vertex of $b$ and $w$ the final vertex of $a$. By Lemma~\ref{lem_renan}, to prove that $\g/\hh$ is a linked net we need only prove that 
$\vp_b(\vp_{ab}^{-1}(W_w))\subseteq W_u$. 

If $\vp_{ab}(\vp_{ab}^{-1}(W_w))=0$ then $\vp_{ab}^{-1}(W_w)=\text{Ker}(\vp_{ab})$, and hence 
$\vp_b(\vp_{ab}^{-1}(W_w))=0$ because $\g$ is a linked net, by Lemma~\ref{lem_renan}. Since 
$W_w$ is simple, we may now assume $\vp_{ab}(\vp_{ab}^{-1}(W_w))=W_w$. In particular, 
$\vp^w_u(W_w)=0$. Since $\hh$ is exact and simple, it follows that $W_w=\vp_a(W_u)$.

By hypothesis, there is a vertex $v$ of $Q$ such that $W_u=\vp^v_u(W_v)$ and $W_v$ is primitive for $\g$. Since $W_w=\vp_a(W_u)$ and $W_w\neq 0$, 
it follows that there is an admissible path $\gamma$ from $v$ to $u$ such that $a\gamma$ is admissible.

If there is an admissible path from $v$ to $z$ passing through $u$, we have
$W_u=\vp_{b}(W_z)$. But then $\vp_{ab}^{-1}(W_w)=W_z+\text{Ker}(\vp_{ab})$, and hence, since $\g$ is a linked net,
$$
\vp_{b}(\vp_{ab}^{-1}(W_w)) = \vp_b(W_z)+\vp_b(\text{Ker}(\vp_{ab}))=W_u+0=W_u.
$$

If not, $\mathfrak a\not\in J$, where $J$ is the essential type of $\ga$. Put $J':=T-J$ and $K:=T-\{\mathfrak a\}$. Since $W_u=\vp^v_u(W_v)$, we must have 
\begin{equation}\label{Wu1}
W_u\subseteq\text{Ker}(\vp^u_{J'}).
\end{equation}
In addition, $\vp_b$ induces the first map in the composition
$$
\vp_{ab}^{-1}(W_w) \longrightarrow \vp_a^{-1}(W_w) \longrightarrow \frac{\vp_a^{-1}(W_w)}{\text{Ker}(\vp_a)},
$$
whose kernel is $\text{Ker}(\vp_{ab})$. Since $W_w=\vp_a(W_u)$, the target of the composition is isomorphic to $W_w$, whence simple. It follows that $\vp_{ab}^{-1}(W_w)=\text{Ker}(\vp_{ab})$ or 
$$
\vp_a^{-1}(W_w)=\text{Ker}(\vp_a)+\vp_b(\vp_{ab}^{-1}(W_w)).
$$
The former cannot occur since $\vp_{ab}(\vp_{ab}^{-1}(W_w))=W_w\neq 0$. Thus the latter holds and yields
\begin{equation}\label{Wu2}
W_u\subseteq \vp_a^{-1}(W_w)=\text{Ker}(\vp_a)+\vp_b(\vp_{ab}^{-1}(W_w))\subseteq 
\text{Ker}(\vp_a)+\text{Ker}(\vp^u_K).
\end{equation}
We conclude from \eqref{Wu1} and \eqref{Wu2} that
$$
W_u\subseteq
\Big(\text{Ker}(\vp^u_K)+\text{Ker}(\vp^u_{\{\mathfrak a\}})\Big)\cap 
\text{Ker}(\vp^u_{J'}).
$$

Since $\g$ satisfies the intersection property at $u$, it follows that 
$$
W_u\subseteq\text{Ker}(\vp^u_{K\cap J'})+\text{Ker}(\vp^u_{\{\mathfrak a\}}), 
$$
and since $\g$ is exact, we get
$$
W_u\subseteq\text{Im}(\vp^{J\cup\{\mathfrak a\}}_u)+\text{Im}(\vp^K_u).
$$
Since $W_u=\vp^v_u(W_v)$, it follows now from Lemma~\ref{SumIma_SumIma} that
$$
W_v\subseteq\text{Im}(\vp^{\{\mathfrak a\}}_v)+\text{Im}(\vp^{J'-\{\mathfrak a\}}_v).
$$
Since $W_v$ is primitive for $\g$, we must have $J'=\{\mathfrak a\}$. But in this case 
$\varphi_a(W_u)=\vp_a\vp^v_u(W_v)=0$, a contradiction with $W_w=\vp_a(W_u)$.
\end{proof}

Recall that an Abelian category is called semisimple if all objects are semisimple, that is, finite direct sums of simple objects. Though the summands are not determined, their number is. The category of finite-dimensional vector spaces is semisimple.

\begin{proposition}\label{exiso} Let $\g$ be a weakly linked net of objects in a $k$-linear Abelian category $\mathcal A$. If $\g$ is exact then all its associated objects have the same class in the Grothendieck group of $\mathcal A$. In particular, if $\mathcal A$ is semisimple, they are isomorphic, $\g$ is pure and every monomorphism between associated objets is an isomorphism.
\end{proposition}

\begin{proof} For each two distinct vertices $u,v$ there is a concatenation of nontrivial simple admissible paths connecting $u$ to $v$. It is thus enough to show that the objects associated to neighboring vertices have the same class in the Grothendieck group of $\mathcal A$. But, given neighboring vertices $u$ and $v$, we have that $V_u$ has the same class as 
$\text{Ker}(\vp^u_v)\oplus \text{Im}(\vp^u_v)$, which is equal to $\text{Ker}(\vp^v_u)\oplus \text{Im}(\vp^v_u)$ by the exactness of $\g$, which has the same class of $V_v$.

If $\mathcal A$ is semisimple then $\mathcal A$ is split, thus objects with the same class in the Grothendieck group are isomorphic. Furthermore, an endomorphism in the category is an epimorphism (resp.~monomorphism) if and only if it is an isomorphism, whence the last statement of the proposition.
\end{proof}

\begin{theorem}\label{thm_geneal_case}
   A locally finite exact linked net over a $\Z^n$-quiver of objects in a $k$-linear semisimple Abelian category is a finite direct sum of simple locally finite exact linked nets if and only if it satisfies the intersection property at every vertex.
\end{theorem}

\begin{proof} Let $\g$ be the linked net of the statement. A direct sum of weakly linked nets satisfies the intersection property at a given vertex, if and only if so does each summand. Thus the ``only if" statement of the theorem follows from Proposition~\ref{SB_inter}, as simple nets are binary.
  
  Conversely, assume $\g$ satisfies the intersection property at every vertex. We may assume $\g$ is nontrivial. We will argue by induction on the number of simple summands of each object associated to $\g$. If all its associated objects are simple, there is nothing to prove.
  
  At any rate, by Lemma~\ref{subrep_prim}, there is a vertex $v$ such that $V^\g_v\neq\sum\text{Im}(\vp^{\mathfrak a}_v)$, where the $\mathfrak a$ run over the arrow types of the quiver. Since the category is semisimple, thus split, there is a simple subobject $W\subseteq V^\g_v$ such that $W\cap(\sum\text{Im}(\vp^{\mathfrak a}_v))=0$. Let $\hh$ be the subnet generated by $W$. Then $\hh$ is faithfully generated by a vertex by Lemma~\ref{subrep_dim_1}, thus simple. By Proposition~\ref{simpleimpliesexact}, the net $\hh$ is locally finite and exact.
  
  The quotient $\g/\hh$ is thus an exact linked net by Proposition~\ref{subrep_quorep}, and satisfies the intersection property at every vertex by Proposition~\ref{thm_quo_inter_pro}. Its associated objects are direct sums of a smaller number of simple objects, when compared to $\g$. By induction, there are 
  simple locally finite exact linked nets $\hh'_1,\dots,\hh'_r$ such that 
  $\g/\hh\cong\hh'_1\oplus\cdots\oplus\hh'_r$.
  
  By Proposition~\ref{simpleimpliesexact}, for each $i=1,\dots,r$ there is a vertex $v_i$ of $Q$ faithfully generating $\hh_i'$. Since the category is split, the quotient map 
  $V^\g_{v_i}\to V^\g_{v_i}/V^{\hh}_{v_i}$ splits, and thus there is a simple suboject 
  $W_i\subseteq V^\g_{v_i}$ isomorphic to $V^{\hh'_i}_{v_i}$ under the quotient. Let $\hh_i$ be the subnet generated by $W_i$. Since $v_i$ faithfully generates $\hh_i'$, and $\g/\hh$ is exact, we must have that $V^{\hh'_i}_{v_i}\cap\text{Im}(\vp^{\g/\hh}_a)=0$ for every arrow $a$ arriving at $v_i$, and thus also $W_i\cap\text{Im}(\vp^{\g}_a)=0$ for every arrow $a$ arriving at $v_i$. It follows from Lemma~\ref{subrep_dim_1} that $\hh_i$ is faithfully generated by $v_i$, whence simple, and thus $\hh_i$ is locally finite and exact by 
  Proposition~\ref{simpleimpliesexact}. Finally, it is clear that $\g=\hh\oplus\hh_1\oplus\cdots\oplus\hh_r$.
\end{proof}

\begin{Cor} A finitely generated exact linked net of finite-dimensional vector spaces is a finite direct sum of exact linked nets of vector spaces of dimension 1 if and only if it satisfies the intersection property at every vertex.
\end{Cor}

\begin{Cor}\label{to_be_stated} A locally finite exact linked net $\g$ over a $\Z^1$-quiver $Q$ of objects in a $k$-linear semisimple Abelian category is a finite direct sum of simple locally finite exact linked nets.
\end{Cor}

\begin{proof} We need only check that $\g$ satisfies the intersection property at 
every vertex, that is, we need to check \eqref{int} for each vertex $v$ of $Q$ and all collections $I_0,\dots,I_m$ of arrow types. But \eqref{int} holds trivially if $I_j$ is the full set of arrow types for some $j$. We may also assume that $I_j\neq\emptyset$ for each $j$. But then there are only two possibilities for each $I_j$ and the property is easily checked by using Lemma~\ref{ker_ker_ker}.
\end{proof}

\begin{Def} Given a $\mathbb Z^n$-quiver $Q$ and a $k$-linear semisimple Abelian category $\mathcal A$, consider the category $\mathcal C_Q(\mathcal A)$ of locally finite exact linked nets over $Q$ of objects in $\mathcal A$, whose morphisms are the representation maps. We call 
$\mathcal C_Q(\mathcal A)$ the category of linked nets over $Q$ with objects in $\mathcal A$. 
\end{Def} 

\begin{theorem}\label{cat} Let $\g$ be a locally finite exact linked net over a $\mathbb Z^n$-quiver $Q$ of objects in a $k$-linear semisimple Abelian category $\mathcal A$. Then $\g$ is simple in $\mathcal C_Q(\mathcal A)$ if and only if $\g$ is simple. And $\g$ is semisimple in $\mathcal C_Q(\mathcal A)$ if and only if $\g$ satisfies the intersection property at every vertex.
\end{theorem} 

\begin{proof} Suppose first that $\g$ is simple. If an exact subnet $\hh\subseteq\g$ is nontrivial, then $V^\hh_v$ is nonzero for some $v$, whence $V^\hh_v$ is nonzero for every $v$ by Proposition~\ref{exiso}. Since $V^\g_v$ is simple for every $v$, it follows that $\hh=\g$.

Conversely, suppose $\g$ is simple in $\mathcal C_Q(\mathcal A)$. Since $\g$ is locally finite, there is a vertex $v$ of $Q$ which is primitive for $\g$ by Lemma~\ref{subrep_prim}. Then there is a simple subnet $\hh$ of $\g$ faithfully generated by $v$ by Lemma~\ref{subrep_dim_1}. The subnet $\hh$ is locally finite and exact by Proposition~\ref{simpleimpliesexact}. Thus $\hh$ is a subobject of $\g$ in $\mathcal C_Q(\mathcal A)$. Since $\g$ is simple in 
$\mathcal C_Q(\mathcal A)$ we must have $\g=\hh$. Thus $\g$ is simple. 

The second statement of the theorem follows from Theorem~\ref{thm_geneal_case} and the first statement.
\end{proof}

Even though all objects in $\mathcal A$ are semisimple, an element in $\mathcal C_Q(\mathcal A)$ is not necessarily semisimple; see Example~\ref{exac_no_SB} below.

%=============================  Example =================================

\section{Applications}\label{applications}

\begin{proposition}\label{simpletriangle} Let $\g$ be an exact linked net over a $\Z^n$-quiver of objects in a $k$-linear semisimple Abelian category. If $\g$ is generated by a 
polygon then $\g$ is a finite direct sum of simple locally finite exact linked nets.
\end{proposition}

\begin{proof} Since $\g$ is finitely generated, $\g$ is locally finite by Proposition~\ref{noseq}. Let $\{w_0,\dots,w_m\}$ be a $(m+1)$-gon generating $\g$. For convenience, put $w_{m+1}:=w_0$. A polygon is its hull by Proposition~\ref{hullpolygon}. So $\g$ is $1$-generated by $\{w_0,\dots,w_m\}$ by Proposition~\ref{simpfin}. By Lemma~\ref{intH} and Theorem~\ref{thm_geneal_case}, it is enough to check that $\g$ satisfies the intersection property at $w_i$ for $i=0,\dots,m$. 

We may assume $m=n$ and that there are arrows $a_0,\dots,a_n$ of 
pairwise distinct types with $a_i$ connecting $w_i$ to $w_{i+1}$ for $i=0,\dots,n$. For each $i=0,\dots,n$, let $J_i$ be the set of all arrow types but $\mathfrak a_i$, the type of $a_i$. Since $\g$ is $1$-generated by $\{w_0,\dots,w_n\}$, for each vertex $w$ of the cone $C_i:=C_{J_i}(w_i)$, the map $\vp^{w_i}_w$ is an isomorphism. The cones $C_0,\dots,C_n$ form a partition of the set of vertices of $Q$. 

Let $I_0,\dots,I_p$ be proper collections of arrow types. We must show that
\begin{equation}\label{intappl}
\bigg(\sum_{\ell=1}^p\text{Ker}(\varphi_{I_\ell}^v)\bigg)\cap\text{Ker}(\varphi_{I_0}^v) = \sum_{\ell=1}^p\text{Ker}(\varphi_{I_\ell\cap I_0}^v)
\end{equation}
holds for $v=w_0,\dots,w_n$. If $p\leq 1$, the equality holds for each vertex $v$ by 
Lemma~\ref{ker_ker_ker}. Apply induction on $p$. Let $d:=\min\{\# I_\ell\}$. If $d=0$ then 
$I_\ell=\emptyset$ for some $\ell$ and thus $\text{Ker}(\varphi_{I_\ell}^v)=\text{Ker}(\varphi_{I_\ell\cap I_0}^v)=0$. Then \eqref{intappl} holds for each vertex $v$ by induction. Apply now induction on $d$.

For each $\ell=1,\dots,p$, we have that $\vp^v_{I_\ell\cap I_0}$ is an monomorphism if $\vp^v_{I_0}$ or $\vp^v_{I_\ell}$ is, as the former is a factor of either of the latter, and hence $\vp^v_{I_\ell\cap I_0}$ is an isomorphism by Proposition~\ref{exiso}. It follows that \eqref{intappl} holds clearly if $\vp^v_{I_0}$ is an isomorphism, whereas it holds by our induction on $p$ if $\vp^v_{I_\ell}$ is an isomorphism for some $\ell>0$. We may thus restrict to collections $I_\ell$ such that $\vp^v_{I_\ell}$ is not an isomorphism for any $\ell=0,\dots,p$.   

Consider \eqref{intappl} for $v:=w_i$ for each $i=0,\dots,n$. So that $\vp^{w_i}_{I_\ell}$ is not an isomorphism for any $\ell$, we must have that $I_\ell\not\subseteq J_i$, and hence $I_\ell$ must contain $\mathfrak a_i$ for every $\ell$. By our induction on $d$, we may assume that 
\begin{equation}\label{i+1}
\bigg(\sum_{\ell=1}^p\text{Ker}(\varphi_{I_\ell\cap J_i}^{w_{i+1}})\bigg)\cap\text{Ker}(\varphi_{I_0\cap J_i}^{w_{i+1}}) = \sum_{\ell=1}^p\text{Ker}(\varphi_{I_\ell\cap I_0\cap J_i}^{w_{i+1}}).
\end{equation}

We prove \eqref{intappl} now for $v:=w_i$. Since $\mathfrak a_i\in I_\ell$ for each $\ell$, it follows that
$$
\vp^{w_i}_{w_{i+1}}(\text{Ker}(\varphi_{I_\ell}^{w_i}))\subseteq \text{Ker}(\varphi_{I_\ell\cap J_i}^{w_{i+1}}).
$$
It follows from \eqref{i+1} that
$$
\vp^{w_i}_{w_{i+1}}\Bigg(\bigg(\sum_{\ell=1}^p\text{Ker}(\varphi_{I_\ell}^{w_i})\bigg)\cap\text{Ker}(\varphi_{I_0}^{w_i})\Bigg)\subseteq \sum_{\ell=1}^p\text{Ker}(\varphi_{I_\ell\cap I_0\cap J_i}^{w_{i+1}}).
$$
Now, since $\g$ is exact, 
$$
\text{Ker}(\varphi_{I_\ell\cap I_0\cap J_i}^{w_{i+1}})\subseteq\text{Ker}(\varphi_{J_i}^{w_{i+1}})=\text{Im}(\vp^{w_i}_{w_{i+1}}),
$$
and hence
$$
\text{Ker}(\varphi_{I_\ell\cap I_0\cap J_i}^{w_{i+1}})=\vp^{w_i}_{w_{i+1}}\big(\text{Ker}(\varphi_{I_\ell\cap I_0}^{w_i})\big)
$$
for each $\ell>0$. It follows that
$$
\bigg(\sum_{\ell=1}^p\text{Ker}(\varphi_{I_\ell}^{w_i})\bigg)\cap\text{Ker}(\varphi_{I_0}^{w_i})\subseteq\text{Ker}(\vp^{w_i}_{w_{i+1}})+\sum_{\ell=1}^p\text{Ker}(\varphi_{I_\ell\cap I_0}^{w_i}).
$$
But, as $I_\ell\cap I_0$ contains $\mathfrak a_i$ for each $\ell$, we must have that 
$\text{Ker}(\vp^{w_i}_{w_{i+1}})\subseteq\text{Ker}(\varphi_{I_\ell\cap I_0}^v)$. Thus we obtained that the left-hand side of \eqref{intappl} is contained in the right-hand side for $v:=w_i$. The reverse inclusion follows from Lemma~\ref{ker_ker_ker}.
\end{proof}

\begin{Exa}\label{exac_no_SB} 
Here we present a finitely generated exact linked net $\g$ of complex vector spaces of dimension 2 over a $\Z^2$-quiver $Q$ which does not satisfy the intersection property at a vertex, and thus is not a direct sum of two linked nets of vector spaces of dimension 1 by Theorem~\ref{thm_geneal_case}. 
The quiver is 
$Q=Q(v;v_0,v_1,v_2)$ with 
$$
v=(2,2,2),\quad v_1=(-2,1,1),\quad v_2=(1,-2,1),\quad v_3=(1,1,-2).
$$
In Figure~2 we describe the representation only over the collection of effective multidegrees, which $1$-generates it.

\begin{figure}[ht!]
\scalebox{0.8}{\begin{minipage}{\textwidth}
\centering
\begin{tikzpicture}[commutative diagrams/every diagram]
  \node  (g)  at (180:5cm)     {$\g:$}; 
  \node (p0) at (0:0cm)       {$\cc^2_{(2,2,2)}$};
  \node (p1) at (90:3.5cm)      {$\cc^2_{(0,3,3)}$};
  \node (p2) at (90+60:3.5cm)   {$\cc^2_{(1,4,1)}$};
  \node (p3) at (90+2*60:3.5cm) {$\cc^2_{(3,3,0)}$};
  \node (p4) at (90+3*60:3.5cm) {$\cc^2_{(4,1,1)}$};
  \node (p5) at (90+4*60:3.5cm) {$\cc^2_{(3,0,3)}$};
  \node (p6) at (90+5*60:3.5cm) {$\cc^2_{(1,1,4)}$};
  \node (p7) at (90+3*60:6cm)    {$\cc^2_{(6,0,0)}$};
  \node (p8) at (90+5*60:7cm) {$\cc^2_{(0,0,6)}$};
  \node (p9) at (90+60:7cm) {$\cc^2_{(0,6,0)}$};
\path[commutative diagrams/.cd, every arrow, every label]  
  (p0) edge node[fill=white,anchor=center, pos=0.5] {$ \setlength{\arraycolsep}{1pt} \begin{bmatrix} 1 & 0 \\[-5pt] 0 & 1 \end{bmatrix}$} (p1)
  (p2) edge node[fill=white,anchor=center, pos=0.5] {$ \setlength{\arraycolsep}{1pt} \begin{bmatrix} 0 & 1 \\[-5pt] 0 & 1 \end{bmatrix}$} (p0)
  (p0) edge node[fill=white,anchor=center, pos=0.5] {$ \setlength{\arraycolsep}{1pt} \begin{bmatrix} 1 & 0 \\[-5pt] 0 & 1 \end{bmatrix}$} (p3)
  (p4) edge node[fill=white,anchor=center, pos=0.5] {$ \setlength{\arraycolsep}{1pt} \begin{bmatrix} 1 & 0 \\[-5pt] 0 & 0 \end{bmatrix}$} (p0) 
  (p0) edge node[fill=white,anchor=center, pos=0.5] {$ \setlength{\arraycolsep}{1pt} \begin{bmatrix} 1 & 0 \\[-5pt] 0 & 1 \end{bmatrix}$} (p5)
  (p6) edge node[fill=white,anchor=center, pos=0.5] {$ \setlength{\arraycolsep}{1pt} \begin{bmatrix} 0 & 0 \\[-5pt] 0 & 1 \end{bmatrix}$} (p0)
  (p8) edge[bend left=15] node[fill=white, anchor=center, pos=0.5] {$0$} (p6)
  (p6) edge[bend left=15] node[fill=white, anchor=center, pos=0.5] {$\cong$} (p8)
  (p9) edge[bend left=15] node[fill=white, anchor=center, pos=0.5] {$0$} (p2)
  (p2) edge[bend left=15] node[fill=white, anchor=center, pos=0.5] {$\cong$} (p9)
  (p7) edge[bend left=15] node[fill=white, anchor=center, pos=0.5] {$0$} (p4)
  (p4) edge[bend left=15] node[fill=white, anchor=center, pos=0.5] {$\cong$} (p7)
  (p1) edge node[fill=white,anchor=center, pos=0.5] {$ \setlength{\arraycolsep}{1pt} \begin{bmatrix} 1 & -1 \\[-5pt] 0 & 0 \end{bmatrix}$} (p2)
  (p1) edge node[fill=white,anchor=center, pos=0.5] {$ \setlength{\arraycolsep}{1pt} \begin{bmatrix} 1 & 0 \\[-5pt] 0 & 0 \end{bmatrix}$} (p6)
  (p5) edge node[fill=white,anchor=center, pos=0.5] {$ \setlength{\arraycolsep}{1pt} \begin{bmatrix} 1 & 0 \\[-5pt] 0 & 0 \end{bmatrix}$} (p6)
  (p5) edge node[fill=white,anchor=center, pos=0.5] {$ \setlength{\arraycolsep}{1pt} \begin{bmatrix} 0 & 0 \\[-5pt] 0 & 1 \end{bmatrix}$} (p4)
  (p3) edge node[fill=white,anchor=center, pos=0.5] {$ \setlength{\arraycolsep}{1pt} \begin{bmatrix} 0 & 0 \\[-5pt] 0 & 1 \end{bmatrix}$} (p4)
  (p3) edge node[fill=white,anchor=center, pos=0.5] {$ \setlength{\arraycolsep}{1pt} \begin{bmatrix} 1 & -1 \\[-5pt] 0 & 0 \end{bmatrix}$} (p2);
\end{tikzpicture}
\end{minipage}}
\label{exa_linked_net_no_SB}
\caption{An exact linked net which is not semisimple.}
\end{figure}
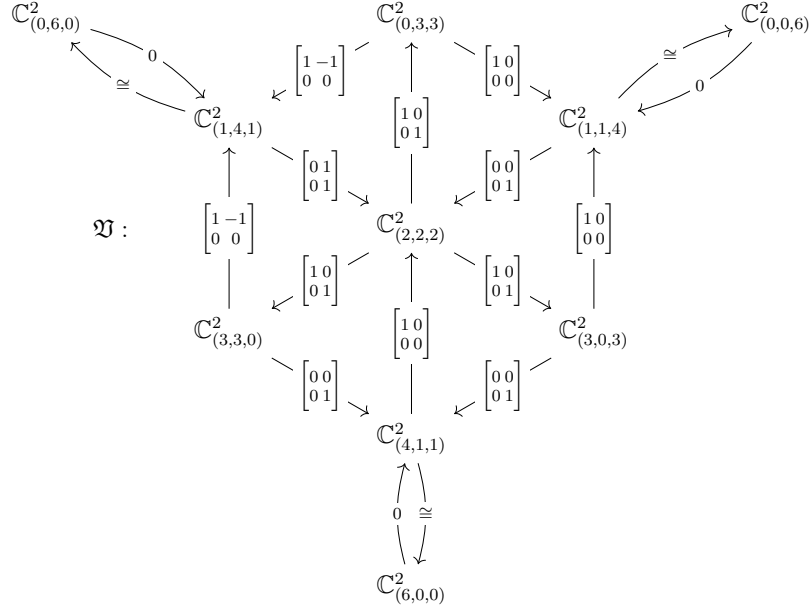
At the vertex $(2,2,2)$ we can verify that 
$$
 \text{Ker}\big(\vp_{(1,4,1)}^{(2,2,2)}\big)=\langle e_1+e_2\rangle, ~ \text{Ker}\big(\vp_{(4,1,1)}^{(2,2,2)}\big)=\langle e_1\rangle ~ \mbox{ and } ~
 \text{Ker}\big(\vp_{(1,1,4)}^{(2,2,2)}\big)=\langle e_2\rangle.
 $$
Since $\text{Ker}\big(\vp^{(2,2,2)}_{(3,3,0)}\big)=0$ and $\text{Ker}\big(\vp^{(2,2,2)}_{(3,0,3)}\big)=0$, we obtain 
 \begin{eqnarray*}
  \Big(\text{Ker}\big(\vp_{(1,4,1)}^{(2,2,2)}\big)+\text{Ker}\big(\vp_{(1,1,4)}^{(2,2,2)}\big)\Big)\cap\text{Ker}\big(\vp_{(4,1,1)}^{(2,2,2)}\big)
  \neq\text{Ker}\big(\vp^{(2,2,2)}_{(3,3,0)}\big)\oplus\text{Ker}\big(\vp^{(2,2,2)}_{(3,0,3)}\big).
 \end{eqnarray*}
 Therefore, $\g$ does not satisfy the intersection property at $(2,2,2)$. 
 
Furthermore, the linked net $\g$ arises from a degeneration of linear series, as we now explain. Let $X$ be the triangle in $\p^2:=\p^2_{\mathbb C}$, given by 
$$
UVW=0.
$$
It has three components, $X_0$, $X_1$ and $X_2$ given by $U=0$, $V=0$ and $W=0$, respectively. Let 
$\fx\subset\p^2\times B$ be the surface defined by $UVW-TF=0$ for a general cubic $F=0$. Here $B:=\text{Spec}(\mathbb C[[T]])$. Since $F$ is general, $\fx$ is regular. Also, the general fiber $\fx_{\eta}$ of the projection $\pi:\fx\to B$ 
is smooth and integral. The special fiber is $X$.

 Consider the invertible sheaf $\mathcal{L}:=\mo_{\fx}(2)$, whose restriction to $X$ has multidegree $(2,2,2)$. The quiver associated to it is the quiver $Q$. 
 
 The coordinates $U,V,W$ of $\p^2$ can be thought as sections of $\mo_{\p^2}(1)$ and restrict to sections of $\mo_{\fx}(1)$ which we denote by $u,v,w$, respectively. Consider the linear system $V_\eta$ of sections of $L_\eta:=\mathcal{L}|_{\fx_{\eta}}$ generated by $u(v+w)$ and $w(v-u)$.
 
For each divisor $D=\sum n_iX_i$ with $\min\{n_i\}=0$, the sheaf $\mathcal L(D)$ may be viewed as a subsheaf of 
$\mo_{\fx}(2+\sum n_i)$, with $B$-flat quotient. Thus, we may think of the limit sections as sections of the larger sheaf $\mo_X(2+\sum n_i)$, which may be represented by polynomials of degree $2+\sum n_i$. In Figure~\ref{degeneration_linear_series} we describe bases for the spaces of limit sections over the collection of effective multidegrees.

\begin{figure}[ht!]
    \scalebox{0.8}{\begin{minipage}{\textwidth}
    \centering
\begin{tikzpicture}[commutative diagrams/every diagram]
 \node  (g)  at (180:5cm)       {$\g:$};
  \node (p0) at (0:0cm)         {$\{u(v+w),w(v-u)\}$};
  \node (p1) at (90:2.8cm)      {$\{u^2(v+w),uw(v-u)\}$};
  \node (p2) at (90+60:2.8cm)   {$\{u^2w(v+w),F(u+w)\}$};
  \node (p3) at (90+2*60:2.8cm) {$\{uw(v+w),w^2(v-u)\}$};
  \node (p4) at (90+3*60:2.8cm) {$\{F(v+w),vw^2(v-u)\}$};
  \node (p5) at (90+4*60:2.8cm) {$\{uv(v+w),vw(v-u)\}$};
  \node (p6) at (90+5*60:2.8cm) {$\{u^2v(v+w),F(v-u)\}$};
  \node (p7) at (90+3*60:5cm) {$\{Fvw(v+w),v^2w^3(v-u)\}$};
  \node (p8) at (90+5*60:5.6cm) {$\{u^3v^2(v+w),Fuv(v-u)\}$};
  \node (p9) at (90+60:5.6cm) {$\{u^3w^2(v+w),Fuw(x+z)\}$};
\path[commutative diagrams/.cd, every arrow, every label]  
%Maps X
(p7) edge node[fill=white,anchor=center, pos=0.5] {$X_0$} (p4)    (p4) edge node[fill=white,anchor=center, pos=0.5] {$X_0$} (p0) 
(p0) edge node[fill=white,anchor=center, pos=0.5] {$X_0$} (p1)    (p5) edge node[fill=white,anchor=center, pos=0.5] {$X_0$} (p6)
(p3) edge node[fill=white,anchor=center, pos=0.5] {$X_0$} (p2)

%Maps Y
(p9) edge node[fill=white,anchor=center, pos=0.5] {$X_1$} (p2)    (p3) edge node[fill=white,anchor=center, pos=0.5] {$X_1$} (p4)  
(p2) edge node[fill=white,anchor=center, pos=0.5] {$X_1$} (p0)    (p0) edge node[fill=white,anchor=center, pos=0.5] {$X_1$} (p5)
(p1) edge node[fill=white,anchor=center, pos=0.5] {$X_1$} (p6)  
%Maps Z    
(p0) edge node[fill=white,anchor=center, pos=0.5] {$X_2$} (p3)    (p6) edge node[fill=white,anchor=center, pos=0.5] {$X_2$} (p0)
(p8) edge node[fill=white,anchor=center, pos=0.5] {$X_2$} (p6)    (p1) edge node[fill=white,anchor=center, pos=0.5] {$X_2$} (p2)  
(p5) edge node[fill=white,anchor=center, pos=0.5] {$X_2$} (p4)  ;
\end{tikzpicture}
\end{minipage}}
\caption{The linked net associated to $\langle u(v+w),w(v-u)\rangle$.}
\label{degeneration_linear_series}
\end{figure}
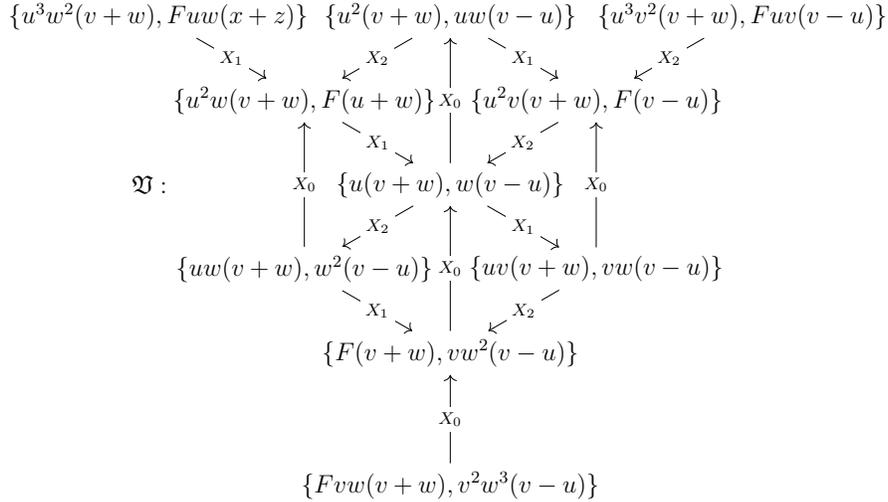

As for the maps between these spaces, the maps labeled $X_0,X_1,X_2$ are multiplication by $u,v,w$, respectively, as long as they map to a space with sections represented by polynomials of higher degree. Otherwise, after multiplication, divide by $F$. Keep in mind that $uvw=0$ on the curve $X$. Using the bases, the maps are represented by the matrices in Figure~2.
\end{Exa}

\bibliographystyle{plain}
\bibliography{Simple_Bases}

\vspace{0.5cm}

{\smallsc Instituto de Matem\'atica Pura e Aplicada, 
Estrada Dona Castorina 110, 22460-320 Rio de Janeiro RJ, Brazil}

{\smallsl E-mail address: \small\verb?esteves@impa.br?}

\vspace{0.2cm}

{\smallsc Instituto de Matem\'atica Pura e Aplicada, 
Estrada Dona Castorina 110, 22460-320 Rio de Janeiro RJ, Brazil}

{\smallsl E-mail address: \small\verb? silvase@impa.br?}

\vspace{0.2cm}

{\smallsc Instituto de Matem\'atica Pura e Aplicada, 
Estrada Dona Castorina 110, 22460-320 Rio de Janeiro RJ, Brazil}

{\smallsl E-mail address: \small\verb?renanmath@ufc.br?}
\end{document}